\documentclass[12pt,a4paper]{article}
\usepackage[utf8]{inputenc}
\usepackage{amsmath, amssymb, amsthm}
\usepackage{geometry}
\usepackage{hyperref}
\usepackage{graphicx}           
\usepackage{booktabs}           
\usepackage{caption}            
\usepackage{subcaption}         
\usepackage{float}              
\usepackage{array}              

\title{A Spline-based Physics-Informed Numerical Scheme: Accurate Smooth Solutions for Differential Equations}
\author{Ayman Mourad$^*$ and Fatima Mrou\'e$^{**}$\\
$^*$Department of Mathematics - Faculty of Sciences (I) - Lebanese University\\
\href{mailto:ayman.mourad@ul.edu.lb}{ayman.mourad@ul.edu.lb}\\
$^{**}$Department of Mathematics and Center of Advanced Mathematical Sciences-\\
American University of Beirut\\
\href{mailto:fm47@aub.edu.lb}{fm47@aub.edu.lb}
}
\date{}

\newcommand{\R}{\mathbb R}

\begin{document}

\maketitle
\begin{abstract}
The rise of Physics-Informed Neural Networks (PINNs) has popularized the concept of solving differential equations via residual minimization. However, neural networks are often viewed as ``black boxes" requiring significant computational overhead and stochastic optimization. Moreover, PINNs typically treat boundary conditions (BCs) as ``soft constraints" within the loss function and this makes the optimization process struggling to enforce the BCs properly. This paper introduces the \textbf{Spline-based Physics-Informed Numerical Scheme (SPINS)}, a numerical framework designed to solve both initial and boundary value problems of ordinary differential equations (ODEs). By replacing the neural network architecture of traditional PINNs with a structured spline basis, SPINS achieves high accuracy and interpretability with a minimal parameter set. In addition, the BCs are automatically satisfied from the choice of the splines architecture. Therefore, SPINS provides smooth numerical solutions for ODEs allowing analytical differentiation. Moreover, SPINS benefits from the automatic differentiation where computing the gradient of the physics-informed loss function is an easy task making the optimization process very fast using gradient-based optimizers such as the L-BFGS-B algorithm. We demonstrate the efficacy of SPINS on nonlinear second order ODEs with several choices of BCs using cubic and quintic interpolating splines and present its natural extension to high order ODEs.
\end{abstract}


\section{Introduction}
For decades, the numerical solution of differential equations has been dominated by mesh-based methods such as the Finite Element Method (FEM), Finite Difference Method (FDM), and Finite Volume Method (FVM), \cite{zienkiewicz2013finite, strang1973analysis, leveque2007finite}. While these techniques are mathematically mature and computationally efficient, they possess inherent limitations in the context of modern ``Physics-Informed" workflows.
Indeed, traditional solvers provide solutions only at specific nodal points or elements. Obtaining values between these points requires interpolation (often low-order), which can introduce artifacts and loss of precision. Moreover, because the solution is discrete, it is difficult to verify if the physics - the differential equation itself - is satisfied throughout the entire domain. Evaluating the residual of the differential equation requires numerical differentiation of the discrete data, which is often noisy and sensitive to grid resolution.\\
On the other hand, the recent shift toward Physics-Informed Neural Networks (PINNs), \cite{raissi2019physics, karniadakis2021physics}, promised a revolution by representing the solution as a continuous, differentiable surrogate model (a neural network). However, as the literature has matured, several critical ``bottlenecks" have been identified. First, PINNs typically treat the boundary conditions (BCs) as ``soft constraints" by adding them to the loss function. This creates a multi-objective optimization problem where the network must balance the differential equation residual against the boundary error. In many stiff or high-frequency problems, the optimizer fails to find a solution that strictly honors the BCs, leading to physically inconsistent results.
In addition, unlike finite elements or finite differences for instance, the weights of a neural network lack geometric meaning. It is difficult to perform local refinement or interpret how specific network parameters influence the local curvature of the solution.
Furthermore, training a PINN requires solving a high-dimensional, non-convex optimization problem, which is significantly slower than the linear or quasi-linear systems found in classical numerical analysis.\\

On the other hand, spline functions represent a powerful alternative to global polynomial interpolation. As established in the foundational work of De Boor \cite{deboor1978practical}, B-splines provide a flexible basis that ensures local control and high-order continuity across a computational domain. In a recent paper \cite{CheaitoMourad}, a recursive algorithm for the construction of interpolating splines of degree $q\geq 2$ has been developped. The algorithm is based on generating $q$ auxiliary splines with prescribed derivatives at the initial knot $x_0$, by using Taylor expansions at knots and the propagation of derivative information, after which any target interpolating spline is obtained as an appropriate linear combination of these auxiliary splines. It is very efficient because we do not need anymore to solve linear systems to compute the splines. Therefore, by construction, the interpolating splines are enforced to satisfy any required conditions. Here we see the tight relationship between BCs of a differential equation and the conditions that we impose on an interpolating spline. As example, Neumann BCs of a second order ODE are intimately linked to the construction of clamped cubic splines where first order derivatives are imposed at both extreme interpolation knots.\\

The concepts of both PINNS and spline interpolation led us to the introduction of the Spline-based Physics-Informed Numerical Scheme (SPINS). SPINS combines the best of both worlds: it utilizes the analytical differentiability of high-order splines (like the $C^2$ cubic and $C^4$  quintic splines) to evaluate the differential equation residuals anywhere in the domain, while maintaining the strict boundary enforcement and computational efficiency of traditional spline interpolation. The spline solution is obtained by minimizing a global residual functional derived from the strong form of the differential operator associated to the ODE. SPINS allows also the computation of the gradient of the loss function thus opening the door to the use of gradient-based optimizers, like the L-BFGS-B algorithm which has shown excellent and fast convergence with the need of only few iterations during the minimization process.

\section{Related Work}
Recent research has increasingly explored the synergy between spline theory and physics-informed frameworks to overcome the limitations of standard coordinate-based neural networks. 
While, the Physics-informed Spline Learning (PiSL) framework by Sun et al. \cite{sun2021physics} demonstrates the robustness of spline-based analytical differentiation for discovering governing equations from noisy, sparse data, 
others like Wandel et al. \cite{wandel2022spline} utilized Hermite spline kernels to bridge grid-based convolutional neural networks representations with continuous physics losses.
Similarly Falini et al. \cite{falini2023splines} have employed Hermite-type quasi-interpolation as a post-processing step to ensure smoothness in PINN-generated domain parameterizations. 
Moreover, the mathematical regularity of splines has been also explored in the work of Cho et al. \cite{cho2021differentiable} which proves that layers of differentiable splines can be seamlessly integrated into gradient-based optimization loops. 
Moving toward more integrated architectures, the Physics-Informed Deep B-Spline Networks (PI-DBSN) framework by Wang et al. \cite{wang2026pidbsn} embeds B-spline functions directly into the learning loop to leverage exact analytical derivatives and satisfy Dirichlet conditions by construction. Unlike earlier spline collocation methods \cite{schiesser2017spline} that relied on rigid grid-based matrix inversions, these modern hybrid approaches treat the spline coefficients as trainable parameters within a ``soft" physics loss, providing a more robust path for solving high-order differential equations and complex dynamical systems.\\

Building on these foundations, the proposed SPINS introduces a novel numerical scheme specifically designed for the robust solution of nonlinear 1D differential equations with possible extensions to multi-dimensional PDEs. While existing literature often uses splines as auxiliary components or neural network layers, SPINS does not consider any network or layers structures but instead it establishes a smooth spline basis as the primary architecture for the solution. By formulating the nonlinear BVP as a direct optimization of the spline coefficients, the proposed scheme ensures high-order regularity and exact satisfaction of boundary conditions. This approach provides a seamless, continuous framework that maintains numerical stability and accuracy even for highly nonlinear operators.
\section{Methodology}
Let $F:[a,b]\times\R^3\longmapsto\R$, $(x,u,v,w)\longmapsto F(x,u,v,w)$, be a regular function and consider a nonlinear second order ODE with general form
\begin{equation}\label{ODE}
F\big(x,u(x),u'(x),u''(x)\big)=0,\quad a\leq x\leq b,
\end{equation}
where $u:[a,b]\longmapsto\R$, $x\longmapsto u(x)$, is the solution obtained after imposing appropriate conditions, i.e. initial conditions (ICs) or BCs.\\

In the sequel of this paper, we will present SPINS in the context of BCs since the case of ICs is much simpler and to avoid redundancy in the derivation of the calculation.\\

The general implicit form of boundary values can be written as follows
\begin{equation}\label{BCs}
g\big(a,u(a),u'(a)\big)=0\quad\textrm{and}\quad h\big(b,u(b),u'(b)\big)=0,
\end{equation}
where $g$ and $h$ are two given functions. For instance, the Dirichlet BCs are of the form: $u(a)=u_a$ and $u(b)=u_b$, and the Neumann BCs can be written as: $u'(a)=v_a$ and $u'(b)=v_b$ where $u_a$, $u_b$, $v_a$ and $v_b$ are given real numbers. The problem formed from (\ref{ODE})-(\ref{BCs}) will be called the boundary value problem (BVP).\\

SPINS approximates the solution $u(x)$ using a spline $S(x)$ of order $q$, for instance using a cubic spline when $q=3$ and a quintic spline when $q=5$. The domain $[a,b]$ is partitioned into $n+1$ knots $\{a=x_0, x_1, \dots, b = x_n\}$. The spline $S(x)$ is uniquely defined by the knot values $\mathbf{y} = [y_0, y_1, \dots, y_n]$ and the BCs at the extreme knots $x_0$ and $x_n$ along with $C^{q-1}$ continuity that is imposed at the internal knots $x_1,\cdots,x_{n-1}$. Here the number of knots need not to be high, it depends on the stiffness of the solution $u$. For instance, $n=5$ is sometimes sufficient to obtain accurate solutions. The number of knots and their locations can be automatically computed, this will be explained in Section \ref{Strategies}.\\

The main idea behind SPINS is to generate first the {\bf necessary} BCs of $S(x)$ at the extreme knots $x_0$ and $x_n$ from the BCs (\ref{BCs}) of the BVP. Therefore, the main unknown becomes the set of knot values $\mathbf{y} = [y_0, y_1, \dots, y_n]$ of the spline $S(x)$. Second, we define a loss function $J(\mathbf{y})$ to be the $L^2$ or $L^{\infty}$ norms of the physics-informed residual as follows
\begin{equation}\label{ObjectiveFct}
J(\mathbf{y}) = \dfrac{1}{2}\|F(\cdot,S,S',S'')\|^2,
\end{equation}
where the derivatives of $S$ are obtained using analytical differentiation. These norms are approximated using a dense set of collocation points $(\xi_j)_{1\leq j\leq N}$ chosen in the interval $[a,b]$. For instance, in practice, the evaluation of $J$ by the $L^2$-norm is done using a discrete sum of squares as follows
\begin{equation}\label{ObjectiveFctL2}
J(\mathbf{y}) \approx \dfrac{1}{2N}\sum_{j=1}^{N}F\Big(\xi_j,S(\xi_j),S'(\xi_j),S''(\xi_j)\Big)^2.
\end{equation}
Consequently, by construction, the spline $S(x)$ already satisfies the BCs (\ref{BCs}) of the BVP. The spline solution of SPINS is obtained from the optimal solution $\mathbf{y}^*$ that we get by solving the minimization problem whose objective function is $J$ defined in (\ref{ObjectiveFct}) or (\ref{ObjectiveFctL2}).\\

The minimization problem is defined on either $\R^{n-1}$, $\R^{n}$ or $\R^{n+1}$ according to which information on $u(a)$ and $u(b)$ are available; both or one of them, or none of them respectively.

\section{Splines Calculation Based on the BVP}\label{Splines}
Let \( (x_0, y_0), \ldots, (x_n, y_n) \) be a set of data points such that $x_0<x_1<\cdots<x_n$. For a given vector $(v_0,v_n)\in\R^2$ there exists a unique cubic spline \( S(x) \) interpolating these points, i.e. $S(x_i) = y_i$ for all $0\leq i\leq n$, such that $S'(x_0) = v_0$ and $S'(x_n)=v_n$. Similarly, for a given vector $(v_0,w_0,v_n,w_n)\in\R^4$ there exists a unique quintic spline \( S(x) \) interpolating these points such that $S'(x_0) = v_0$, $S''(x_0)=w_0$, $S'(x_n) = v_n$ and $S''(x_n)=w_n$. These are called the clamped cubic and quintic splines.\\

We denote by $W(v_0,v_n)$ the clamped cubic spline interpolating the given knots that verifies $S'(x_0) = v_0$ and $S'(x_n)=v_n$. Similarly, we denote by $W(v_0,w_0,v_n,w_n)$ the clamped quintic spline interpolating the given knots that verifies $S'(x_0) = v_0$, $S''(x_0)=w_0$, $S'(x_n) = v_n$ and $S''(x_n)=w_n$.\\

Moreover, let $W_0=W(0,0)$, $W_1=W(1,0)$, and $W_2=W(0,1)$ for the clamped cubic splines, and $W_0=W(0,0,0,0)$, $W_1=W(1,0,0,0)$, $W_2=W(0,1,0,0)$, $W_3=W(0,0,1,0)$, and $W_4=W(0,0,0,1)$ for the clamped quintic splines. The splines $W_0,\cdots,W_4$ are called {\bf the basic clamped spline functions} of the given data points.\\

Therefore, by the definition of the $W_i$'s, we have the following:
\begin{itemize}
\item Cubic spline: 
$$\left\{
\begin{array}{l}
W_0'(a) = 0, W_1'(a) = 1, W_2'(a) = 0\medskip\\
W_0'(b) = 0, W_1'(b) = 0, W_2'(b) = 1.
\end{array}\right.$$
\item Quintic spline: 
$$\left\{
\begin{array}{l}
W_0'(a) = 0, W_1'(a) = 1, W_2'(a) = 0, W_3'(a) = 0, W_4'(a) = 0\medskip\\
W_0''(a) = 0, W_1''(a) = 0, W_2''(a) = 1, W_3''(a) = 0, W_4''(a) = 0\bigskip\\
W_0'(b) = 0, W_1'(b) = 0, W_2'(b) = 0, W_3'(b) = 1, W_4'(b) = 0\medskip\\
W_0''(b) = 0, W_1''(b) = 0, W_2''(b) = 0, W_3''(b) = 0, W_4''(b) = 1.
\end{array}\right.$$
\end{itemize}
It is easy to show that, see \cite{CheaitoMourad}, the clamped splines can be expressed as a combination of the basic clamped spline functions as follows
$$W(v_0,v_n) = (1-v_0-v_n)W_0+v_0W_1+v_nW_2,$$
and
$$W(v_0,w_0,v_n,w_n) = (1-v_0-w_0-v_n-w_n)W_0+v_0W_1+w_0W_2+v_nW_3+w_nW_4.$$
Furthermore, any other spline $S$ (cubic or quintic) that interpolates these data points can be uniquely written as a combination of the basic clamped spline functions as follows
$$S = (1-\alpha_1-\alpha_2)W_0+\alpha_1W_1+\alpha_2W_2,$$
where 
$$S'(a) = \alpha_1,\quad S'(b) = \alpha_2,$$
or
$$S = (1-\alpha_1-\alpha_2-\alpha_3-\alpha_4)W_0+\alpha_1W_1+\alpha_2W_2+\alpha_3W_3+\alpha_4W_4,$$
where 
$$S'(a) = \alpha_1,\quad S''(a) = \alpha_2,\quad S'(b) = \alpha_3,\quad S''(b) = \alpha_4.$$

In summary, if we have as given the data points, then we can first compute the basic clamped spline functions, and then any other spline interpolating these points can be uniquely determined if we know its coefficients $\alpha_i$'s. And here it comes the role of the BCs (\ref{BCs}) of the BVP to help in finding these coefficients.\\

Now we will consider several case studies of the BVP in order to compute both the corresponding cubic and quintic splines.
\subsection{Cubic spline with Neumann BCs}
This is the simplest scenario where the original BVP consists of imposing the first order derivative at $a$ and $b$ simultaneously. In other words, we do have the values of $u'(a)$ and $u'(b)$. So here the corresponding cubic spline is the clamped spline $W(u'(a),u'(b))$, i.e.
$$S = (1-u'(a)-u'(b))W_0+u'(a)W_1+u'(b)W_2,$$
where the clamped basic spline functions interpolate the data points $(x_i,y_i)$, $i=0,\cdots,n$.\\
The optimization problem is defined over $\R^{n+1}$ for the unknown vector $\mathbf{y} = [y_0,\cdots,y_n]$.

\subsection{Cubic spline with Dirichlet BCs}
Now the BCs consist of the values of $u$ at the boundaries, i.e. $u(a)$ and $u(b)$ are given. So the optimization problem is defined over $\R^{n-1}$ and the unknown vector is $\mathbf{y} = [y_1,\cdots,y_{n-1}]$ since $y_0 = u(a)$ and $y_n = u(b)$.\\
The corresponding spline interpolates the data points $(x_i,y_i)$, $i=0,\cdots,n$ and it is expressed as follows
$$S = (1-\alpha_1-\alpha_2)W_0+\alpha_1W_1+\alpha_2W_2$$
Plugging $S$ into the ODE (\ref{ODE}) for $x=a$ and $x=b$, we get the following two equations:
\begin{equation}\label{CS_Diri}
F(a,u(a),S'(a),S''(a)) = 0\quad\textrm{and}\quad F(b,u(b),S'(b),S''(b)) = 0,
\end{equation}
where 
$$S'(a) = \alpha_1,\quad S'(b) = \alpha_2,$$
$$S''(a) = (1-\alpha_1-\alpha_2)W_0''(a)+\alpha_1W_1''(a)+\alpha_2W_2''(a),$$
and 
$$S''(b) = (1-\alpha_1-\alpha_2)W_0''(b)+\alpha_1W_1''(b)+\alpha_2W_2''(b).$$
Therefore we have a nonlinear system of two equations with two unknowns $\alpha_1$ and $\alpha_2$ that can be easily solved using for instance Newton-Raphson's method.

\subsection{Quintic spline with Neumann BCs}\label{QuinticNeumann}
In this case we have the values of $u'(a)$ and $u'(b)$. Therefore the optimization problem is defined over $\R^{n+1}$ for the unknown vector $\mathbf{y} = [y_0,\cdots,y_n]$.\\

However, in order to compute a clamped quintic spline we need the information of $u''(a)$ and $u''(b)$ which are not explicitly available.
So here the corresponding quintic spline is given by
$$S = (1-\alpha_1-\alpha_2-\alpha_3-\alpha_4)W_0+\alpha_1W_1+\alpha_2W_2+\alpha_3W_3+\alpha_4W_4$$
where $\alpha_1 = u'(a)$ and $\alpha_3 = u'(b)$ are given.\\
In order to compute $\alpha_2$ and $\alpha_4$, we plug $S$ into the ODE (\ref{ODE}) and get the following system
\begin{equation}\label{QS_Neum}
F(a,S(a),u'(a),S''(a)) = 0\quad\textrm{and}\quad F(b,S(b),u'(b),S''(b)) = 0,
\end{equation}
where 
$$S(a) = y_0,\quad S(b) = y_n,\quad S''(a) = \alpha_2,\quad S''(b) = \alpha_4.$$
This leads to the following equations:
\begin{equation}\label{QS_Neum_Simple}
F(a,y_0,u'(a),\alpha_2) = 0\quad\textrm{and}\quad F(b,y_n,u'(b),\alpha_4) = 0.
\end{equation}
Thus we get two independent nonlinear equations with unknowns $\alpha_2$ and $\alpha_4$. These equations can be easily solved numerically in the case of a fully nonlinear ODE using for instance Newton-Raphson's method. But in the case of a linear, semi-linear or quasi-linear ODEs, $\alpha_2$ can be expressed explicitly in terms of $a$, $y_0$ and $u'(a)$, and $\alpha_4$ in terms of $b$, $y_n$ and $u'(b)$.
\subsection{Quintic spline with Dirichlet BCs}
Now we know the values of $u(a)$ and $u(b)$. So the optimization problem is defined over $\R^{n-1}$ and the unknown vector is $\mathbf{y} = [y_1,\cdots,y_{n-1}]$ with $y_0 = u(a)$ and $y_n = u(b)$.\\
The corresponding quintic spline interpolates the data points $(x_i,y_i)$, $i=0,\cdots,n$ and it is expressed as follows
$$S = (1-\alpha_1-\alpha_2-\alpha_3-\alpha_4)W_0+\alpha_1W_1+\alpha_2W_2+\alpha_3W_3+\alpha_4W_4$$
The four parameters $\alpha_1$, $\alpha_2$, $\alpha_3$ and $\alpha_4$ are all unknowns.
Plugging $S$ into the ODE (\ref{ODE}) for $x=a$ and $x=b$, provides only the following two equations:
\begin{equation}\label{QS_Diri}
F(a,u(a),S'(a),S''(a)) = 0\quad\textrm{and}\quad F(b,u(b),S'(b),S''(b)) = 0,
\end{equation}
where 
$$S'(a) = \alpha_1,\quad S''(a) = \alpha_2,\quad S'(b) = \alpha_3,\quad S''(b) = \alpha_4.$$
Thus the equations in (\ref{QS_Diri}) reduce to the following
\begin{equation}\label{QS_Diri_Simple}
F(a,u(a),\alpha_1,\alpha_2) = 0\quad\textrm{and}\quad F(b,u(b),\alpha_3,\alpha_4) = 0.
\end{equation}
In order to determine these four unknowns, we still have to consider two more equations. To this end, and since our solution is assumed to be smooth, we differentiate our ODE (\ref{ODE}) with respect to $x$ and get the following:
\begin{equation}\label{d_ODE_dx}
\dfrac{d}{dx}
F\big(x,u(x),u'(x),u''(x)\big) = \dfrac{\partial F}{\partial x}+\dfrac{\partial F}{\partial u}u'(x)+\dfrac{\partial F}{\partial v}u''(x)+\dfrac{\partial F}{\partial w}u'''(x)=0.
\end{equation}
Now plugging $S$ into (\ref{d_ODE_dx}) for $x=a$ and $x=b$ respectively, we obtain
\begin{equation}\label{QS_Diri_Augmented}
\begin{array}{l}
\dfrac{\partial F}{\partial x}\big(a,y_0,\alpha_1,\alpha_2\big)+\alpha_1\dfrac{\partial F}{\partial u}\big(a,y_0,\alpha_1,\alpha_2\big)+\alpha_2\dfrac{\partial F}{\partial v}\big(a,y_0,\alpha_1,\alpha_2\big)
+S'''(a)\dfrac{\partial F}{\partial w}\big(a,y_0,\alpha_1,\alpha_2\big)=0,\medskip\\
\textrm{and}\bigskip\\
\dfrac{\partial F}{\partial x}\big(b,y_n,\alpha_3,\alpha_4\big)+\alpha_3\dfrac{\partial F}{\partial u}\big(b,y_n,\alpha_3,\alpha_4\big)+\alpha_4\dfrac{\partial F}{\partial v}\big(b,y_n,\alpha_3,\alpha_4\big)
+S'''(b)\dfrac{\partial F}{\partial w}\big(b,y_n,\alpha_3,\alpha_4\big)=0,
\end{array}
\end{equation}
where
$$S'''(a) = (1-\alpha_1-\alpha_2-\alpha_3-\alpha_4)W_0'''(a)+\alpha_1W_1'''(a)+\alpha_2W_2'''(a)+\alpha_3W_3'''(a)+\alpha_4W_4'''(a),$$
and 
$$S'''(b) = (1-\alpha_1-\alpha_2-\alpha_3-\alpha_4)W_0'''(b)+\alpha_1W_1'''(b)+\alpha_2W_2'''(b)+\alpha_3W_3'''(b)+\alpha_4W_4'''(b).$$

Therefore we have an augmented nonlinear system (\ref{QS_Diri_Simple})-(\ref{QS_Diri_Augmented}) of four equations with four unknowns $\alpha_1$, $\alpha_2$, $\alpha_3$ and $\alpha_4$ that can be solved using for instance Newton-Raphson's method.

\subsection{Arbitrary BCs}
Suppose that the BCs are given in the general form (\ref{BCs}), where both $u$ and $u'$ appear explicitly in the expressions of $g$ and $h$ as in the case of Robin's BCs which are given as $u'(a)+\gamma_a u(a) = c_a$ and $u'(b)+\gamma_b u(b) = c_b$ where $\gamma_a$,$\gamma_b$, $c_a$ and $c_b$ are real constants. In this case, the loss function $J$ will be defined on $\R^{n+1}$ where the unknown vector is $\mathbf{y} = [y_0,\cdots,y_{n}]$.\\

During the optimization process, and through the evaluation of $J(\mathbf{y})$, we solve simultaneously, for the unknowns $\alpha_1$ and $\alpha_2$ in terms of $y_0$ and $y_n$ respectively, the following nonlinear equations
\begin{equation}\label{BCs_1-2}
g\big(a,y_0,\alpha_1\big)=0\quad\textrm{and}\quad h\big(b,y_n,\alpha_2\big)=0.
\end{equation}
If we are considering cubic splines, then we have all what we need to compute the corresponding clamped cubic spline $S = W(\alpha_1,\alpha_2)$ and its first and second order derivatives.\\

Otherwise, if we are considering quintic splines, then in addition to $\alpha_1$ and $\alpha_2$ which are obtained from (\ref{BCs_1-2}), we need to solve simultaneously, for the unknown $\alpha_3$ in terms of $a$, $y_0$ and $\alpha_1$, and for the unknown $\alpha_4$ in terms of $b$, $y_n$ and $\alpha_2$ respectively, the following two additional nonlinear equations 
\begin{equation}\label{BCs_3-4}
F\big(a,y_0,\alpha_1,\alpha_3\big) = 0,\quad\textrm{and}\quad F\big(b,y_n,\alpha_2,\alpha_4\big) = 0.
\end{equation}
Consequently, we have all of the $\alpha_i$'s, so we can compute the corresponding clamped quintic spline $S = W(\alpha_1,\alpha_2,\alpha_3,\alpha_4)$ and its derivatives up to order $4$.\\

It is worth to mention that in the case where the BCs $u'(a)$ and $u'(b)$ can be expressed linearly or nonlinearly but explicitly in terms of $u(a)$ and $u(b)$ respectively, then the evaluation of $\alpha_1$ and $\alpha_2$ becomes explicit and there is no need to solve nonlinear equations. This is for instance the case of Robin's BCs. Similarly, when the ODE (\ref{ODE}) is linear, semi-linear or quasi-linear, then the evaluation of $\alpha_3$ and $\alpha_4$ becomes explicit.\\

If the BCs are not given at $a$ and $b$ separately and they are mixed conditions for both boundary points $a$ and $b$, i.e. a boundary condition that involves $u$ and its derivatives at both ends $a$ and $b$ in the same equation, then instead of System (\ref{BCs_1-2})-(\ref{BCs_3-4}), we will be dealing with a coupled system of four equations with the four unknowns the $\alpha_i$'s.

\section{High Order BVPs}
In the previous sections, we presented the evaluation of splines based on BCs in the context of a second order nonlinear BVP. In this section we will show how this can be easily extended to high order BVPs. Indeed, the use of quintic splines is applicable to BVPs up to order $4$.
In particular, we will present the derivation of the spline calculation for the case of a fourth order BVPs.\\

Let $F:[a,b]\times\R^5\longmapsto\R$, $(x,u,v,w,p,q)\longmapsto F(x,u,v,w,p,q)$, be a regular function and consider a nonlinear fourth order ODE with general form
\begin{equation}\label{ODE_high}
F\big(x,u(x),u'(x),u''(x),u'''(x),u^{(4)}(x)\big)=0,\quad a\leq x\leq b,
\end{equation}
completed with boundary values written in a general implicit form as follows
\begin{equation}\label{BCs_4th}
\left\{\begin{array}{l}
g_1\big(a,u(a),u'(a),u''(a)\big)=0,\quad g_2\big(a,u(a),u'(a),u''(a)\big)=0,\medskip\\
h_1\big(b,u(b),u'(b),u''(b)\big)=0,\quad h_2\big(b,u(b),u'(b),u''(b)\big)=0,
\end{array}\right.
\end{equation}
where $g_1$, $g_2$, $h_1$ and $h_2$ are given functions.\\

If we denote by $S$ the quintic spline solution provided by SPINS, then we have
$$S = (1-\alpha_1-\alpha_2-\alpha_3-\alpha_4)W_0+\alpha_1W_1+\alpha_2W_2+\alpha_3W_3+\alpha_4W_4,$$
where $S'(a) = \alpha_1$, $S''(a) = \alpha_2$, $S'(b) = \alpha_3$, and $S''(b) = \alpha_4$.\\
Plugging $S$ into the BCs (\ref{BCs_4th}), we obtain four equations with unknowns the coefficients $\alpha_i$'s
\begin{equation}\label{BCs_4th_SPINS}
\begin{array}{ll}
\left\{\begin{array}{l}
g_1\big(a,y_0,\alpha_1,\alpha_2\big)=0,\medskip\\
g_2\big(a,y_0,\alpha_1,\alpha_2\big)=0,
\end{array}\right. &
\left\{\begin{array}{l}
h_1\big(b,y_n,\alpha_3,\alpha_4\big)=0,\medskip\\
h_2\big(b,y_n,\alpha_3,\alpha_4\big)=0.
\end{array}\right.
\end{array}
\end{equation}
System (\ref{BCs_4th_SPINS}) consists of two sets of two equations each, with unknowns $\alpha_1$ and $\alpha_2$ for one set, and $\alpha_3$ and $\alpha_4$ for the other one. The numerical solving of these sets in the worst case of fully nonlinear functions, can be done using Newton-Raphson's method. Otherwise, exact solving is possible when the BCs are not fully nonlinear which is the case in most realistic $4^{th}$ order BVPs. On the other hand, if the BCs contain mixed conditions for both boundary points $a$ and $b$, then the equivalent system similar to (\ref{BCs_4th_SPINS}) will be solved as a whole system of four equations with the four unknowns the $\alpha_i$'s.\\

Notice that for fourth order BVPs, the BCs provide exactly four equations at the boundaries, then we do not need to generate extra equations by differentiating the ODE as we did in the case of second order BVP with Dirichlet BCs and dealing with quintic spline solutions.\\

However, if the BVP is of $3^{rd}$ order, then we have only three equations for the BCs. This means that we need one more equation to solve for the unknowns the $\alpha_i$'s. Such an equation can be obtained by differentiating the ODE and plugging the quintic spline into it at one appropriate endpoint $a$ or $b$. Since the BVP is of the third order, this extra equation will contain a term with the fourth order derivative of $S$, and this is allowed because the quintic spline is of class $C^4$.\\

If the BVP is of order greater than or equal to $5$, then splines of order $q\geq 6$ are required. In this case, we proceed as we did above. The computation of splines of order $q\geq 6$ can be done iteratively based on the ``Forward and Update" algorithm presented in \cite{CheaitoMourad}.
\section{Automatic Differentiation: Towards Gradient-Based Optimization}
In this section, we present the calculation of the gradient of the loss function $J$ using automatic differentiation in a way similar to the backpropagation step used in neural networks.\\
In order to avoid complex notations in the calculation of the gradient, we will do the derivation in the case of cubic splines.\\
Given knots \( x_0 < x_1 < \dots < x_n \) and knot values $\mathbf{y} = [y_0,\cdots,y_n]$, and as mentioned in section \ref{Splines}, the clamped spline interpolating the points $(x_i,y_i)_i$'s with boundary derivatives $v_0$ and $v_n$ is denoted $W(v_0,v_n)$. We denote by $R(v_0,w_0)$ the cubic spline that interpolates the points $(x_i,y_i)_i$'s with initial conditions, the first and second order derivatives at the initial knot $x_0$, equal to $v_0$ and $w_0$ respectively. This spline will be called the {\bf recursive cubic spline} with initial conditions $v_0$ and $w_0$. The restriction of $R(v_0,w_0)$ on each subinterval \( [x_i, x_{i+1}] \) is a cubic polynomial of the form:
\[
S_i(x) = a_i + b_i(x - x_i) + c_i(x - x_i)^2 + d_i(x - x_i)^3
\]
with the following requirements
\begin{itemize}
    \item Interpolation: \( S_i(x_i) = y_i, \quad S_i(x_{i+1}) = y_{i+1} \),
    \item \( C^1 \) continuity: \( S_i'(x_{i+1}) = S_{i+1}'(x_{i+1}) \),
    \item \( C^2 \) continuity: \( S_i''(x_{i+1}) = S_{i+1}''(x_{i+1}) \),
    \item Initial conditions: $S_0'(x_0) = v_0$ and $S_0''(x_0)=w_0$.
\end{itemize}
Basic calculation leads to the following recursion formula, for $1\leq i\leq n-1$,
$$\left\{
\begin{array}{l}
a_i = y_i\medskip\\
b_i = b_{i-1}+2c_{i-1}h_{i-1}+3d_{i-1}h_{i-1}^2\medskip\\
c_i = c_{i-1}+3d_{i-1}h_{i-1}\medskip\\
d_i = \dfrac{1}{h_i^3}\left(y_{i+1}-y_i-b_{i}h_i-c_{i}h_{i}^2\right)
\end{array}\right.$$
where $h_i = x_{i+1}-x_i$, and for $i=0$ we have
$$a_0 = y_0,\quad b_0 = v_0,\quad c_0 = \dfrac{w_0}{2},\quad d_0 = \dfrac{y_1-y_0}{h_0^3}-\dfrac{b_0+c_0h_0}{h_0^2}.$$
For $0\leq i\leq n-1$, we define the vector $\mathbf{V}_i=\left(\begin{array}{c}
a_i\\ b_i\\ c_i\\ d_i
\end{array}\right)$ containing the coefficients of $R(v_0,w_0)$. We have then
$$\mathbf{V}_0=\left(\begin{array}{c}
a_0\\ b_0\\ c_0\\ d_0
\end{array}\right)
=\left(\begin{array}{c}
y_0\\ 0\\ 0\\ \dfrac{y_1-y_0}{h_0^3}
\end{array}\right)
+\left(\begin{array}{c}
0\\ v_0\\ \dfrac{w_0}{2}\\ -\dfrac{2v_0+w_0h_0}{2h_0^2}
\end{array}\right) = \mathbf{A}_0\mathbf{y}+\mathbf{C}_0(v_0,w_0)$$
where $\mathbf{A}_0 = \left(\begin{array}{ccccc}
1 & 0 & 0 & \cdots & 0\\
0 & 0 & 0 & \cdots & 0\\
0 & 0 & 0 & \cdots & 0\\
-\frac{1}{h_0^3} & \frac{1}{h_0^3} & 0 & \cdots & 0
\end{array}\right)$ is a $4\times (n+1)$ matrix.\\

Therefore, we get the following relation
$$\begin{array}{lcl}
\mathbf{V}_1=\left(\begin{array}{c}
a_1\\ b_1\\ c_1\\ d_1
\end{array}\right)
& = &\left(\begin{array}{cccccc}
0 & 1 & 0 & 0 & \cdots & 0\\
0 & 0 & 0 & 0 & \cdots & 0\\
0 & 0 & 0 & 0 & \cdots & 0\\
0 & -\frac{1}{h_1^3} & \frac{1}{h_1^3} & 0 & \cdots & 0
\end{array}\right)
\left(\begin{array}{c}
y_0\\ \vdots\\ \vdots\\ y_n
\end{array}\right)
\medskip\\
& & 
+\left(\begin{array}{cccc}
0 & 0 & 0 & 0\\
0 & 1 & 2h_0 & 3h_0^2\\
0 & 0 & 1 & 3h_0\\
0 &  -\frac{1}{h_1^2} &  -\frac{2h_0+h_1}{h_1^2} & -\frac{3h_0(h_0+h_1)}{h_1^2}
\end{array}\right)\left(\begin{array}{c}
a_0\\ b_0\\ c_0\\ d_0
\end{array}\right)\medskip\\
& = &\mathbf{A}_{1}\mathbf{y}+\mathbf{B}_{0}\mathbf{V}_{0}
\end{array}$$
In other words, we have the following
\begin{equation}\label{RecursiveLinear}
\left\{\begin{array}{l}
\mathbf{V}_0 = \mathbf{A}_0\mathbf{y}+\mathbf{C}_0(v_0,w_0)\medskip\\
\mathbf{V}_{i} = \mathbf{A}_{i}\mathbf{y}+\mathbf{B}_{i-1}\mathbf{V}_{i-1},\,\forall 1\leq i\leq n-1
\end{array}\right.
\end{equation}
where $\mathbf{A}_{i}$ is a $4\times (n+1)$ matrix that has all entries are zeros except $\mathbf{A}_i(1,i+1) = 1$, $\mathbf{A}_i(4,i+1) = -\frac{1}{h_i^2}$ and $\mathbf{A}_i(4,i+2) = \frac{1}{h_i^2}$, and $\mathbf{B}_{i-1}$ is a $4\times 4$ matrix given by 
$$\mathbf{B}_{i-1} = \left(\begin{array}{cccc}
0 & 0 & 0 & 0\\
0 & 1 & 2h_{i-1} & 3h_{i-1}^2\\
0 & 0 & 1 & 3h_{i-1}\\
0 &  -\frac{1}{h_i^2} &  -\frac{2h_{i-1}+h_i}{h_i^2} & -\frac{3h_{i-1}(h_{i-1}+h_i)}{h_i^2}
\end{array}\right).$$

This means that, according to the iterative formula given in (\ref{RecursiveLinear}), the coefficients of the recursive cubic spline $R(v_0,w0)$  depend linearly (up to an additive constant vector) on the knot values vector $\mathbf{y}$. Consequently, the gradients of the vectors $\mathbf{V}_{i}$'s are given by
\begin{equation}\label{RecursiveGradient}
\left\{\begin{array}{l}
\nabla_{\mathbf{y}}\mathbf{V}_0 = \mathbf{A}_0\medskip\\
\nabla_{\mathbf{y}}\mathbf{V}_{i} = \mathbf{A}_{i}+\mathbf{B}_{i-1}\nabla_{\mathbf{y}}\mathbf{V}_{i-1},\,\forall 1\leq i\leq n-1
\end{array}\right.
\end{equation}
It is very important to notice that the $\nabla_{\mathbf{y}}\mathbf{V}_{i}$'s are independent of $v_0$ and $w_0$, and clearly independent of $\mathbf{y}$. They only depend on the $h_i$'s. This means that they will be computed only once when we have the knots positions $x_i$'s.\\

As for the clamped cubic splines, let $R_0 = R(v_0,0)$ and $R_1 = R(v_0,1)$. Using the result from \cite{CheaitoMourad}, we have the following relationship
\begin{equation}\label{ClampedRecursive}
W(v_0,v_n) = \alpha R_0+(1-\alpha)R_1,\,\,\textrm{where}\,\, \alpha = \dfrac{v_n - R_1'(x_n)}{R_0'(x_n) - R_1'(x_n)}. 
\end{equation}
From (\ref{RecursiveGradient}), we can compute $\nabla_{\mathbf{y}}R(v_0,w0)$, $\nabla_{\mathbf{y}}R'(v_0,w0)$ and $\nabla_{\mathbf{y}}R''(v_0,w0)$. Thus $\nabla_{\mathbf{y}}\alpha$ can be evaluated and consequently, $\nabla_{\mathbf{y}}W(v_0,v_n)$ can be obtained from (\ref{ClampedRecursive}). In particular, the gradients of the basic clamped spline functions $W_0$, $W_1$ and $W_2$ and their first and second order derivatives can be computed.\\

Now, considering the loss function $J$, according to (\ref{ObjectiveFctL2}), we write
$$J(\mathbf{y}) = \dfrac{1}{N}\sum_{k=1}^{N}J_k(\mathbf{y})$$
where $J_k(\mathbf{y}) = \dfrac{1}{2} F\big(\xi_k,S(\xi_k),S'(\xi_k),S''(\xi_k)\big)^2$, and $S$ is the unique cubic spline that interpolates the points $(x_i,y_i)_i$'s and verifying the required BCs.\\

Therefore we have
$$\nabla_{\mathbf{y}}J_k(\mathbf{y}) = F\big(\xi_k,S(\xi_k),S'(\xi_k),S''(\xi_k)\big)\nabla_{\mathbf{y}}F\big(\xi_k,S(\xi_k),S'(\xi_k),S''(\xi_k)\big)$$
and
$$\begin{array}{lccl}
\nabla_{\mathbf{y}}F\big(\xi_k,S(\xi_k),S'(\xi_k),S''(\xi_k)\big) & = & &\dfrac{\partial F}{\partial u}\big(\xi_k,S(\xi_k),S'(\xi_k),S''(\xi_k)\big)\nabla_{\mathbf{y}}S(\xi_k)\medskip\\
& & + &\dfrac{\partial F}{\partial v}\big(\xi_k,S(\xi_k),S'(\xi_k),S''(\xi_k)\big)\nabla_{\mathbf{y}}S'(\xi_k)\medskip\\
& & + &\dfrac{\partial F}{\partial w}\big(\xi_k,S(\xi_k),S'(\xi_k),S''(\xi_k)\big)\nabla_{\mathbf{y}}S''(\xi_k).
\end{array}
$$

In order to complete our calculation, it remains to determine the gradients of $S$, $S'$ and $S''$. This can be accomplished by writing $S$ as a combination of the basic clamped spline functions as follows
$$S = (1-\alpha_1-\alpha_2)W_0+\alpha_1W_1+\alpha_2W_2,$$
where 
$$S'(a) = \alpha_1,\quad S'(b) = \alpha_2.$$

In the case of Neumann BCs, the values of $S'(a)$ and $S'(b)$ are given, thus $\alpha_1$ and $\alpha_2$ are given constants independent of $\mathbf{y}$, hence $\nabla_{\mathbf{y}}S = (1-\alpha_1-\alpha_2)\nabla_{\mathbf{y}}W_0+\alpha_1\nabla_{\mathbf{y}}W_1+\alpha_2\nabla_{\mathbf{y}}W_2$. We have similar expressions for $\nabla_{\mathbf{y}}S'$ and $\nabla_{\mathbf{y}}S''$.\\

In the case of Dirichlet BCs, $u(a)$ and $u(b)$ are given. Then the unknown vector becomes $\mathbf{y} = [y_1,\cdots,y_{n-1}]\in\R^{n-1}$. So the matrices $\mathbf{A}_{i}$ are of dimension $4\times(n-1)$ and can be obtained by removing the first and last columns from their definition given above. We should also make an appropriate change to $\mathbf{C}_0(v_0,w_0)$ to account for $y_0 = u(a)$ and add a vector containing $y_n = u(b)$ to $\mathbf{V}_{n-1}$. However, and most importantly, now $\alpha_1$ and $\alpha_2$ satisfy System (\ref{CS_Diri}), thus they depend implicitly on the vector $\mathbf{y}$ and we must compute their gradients $\nabla_{\mathbf{y}}\alpha_1$ and $\nabla_{\mathbf{y}}\alpha_2$. This can be achieved by differentiating the equations in (\ref{CS_Diri}) as follows:
\begin{equation}\label{CS_Diri_Gradient}
\begin{array}{l}
\dfrac{\partial F}{\partial v}\big(a,u(a),\alpha_1,S''(a)\big)\nabla_{\mathbf{y}}\alpha_1
+ \dfrac{\partial F}{\partial w}\big(a,u(a),\alpha_1,S''(a)\big)\nabla_{\mathbf{y}}S''(a) = 0,\medskip\\
\dfrac{\partial F}{\partial v}\big(b,u(b),\alpha_2,S''(b)\big)\nabla_{\mathbf{y}}\alpha_2
+ \dfrac{\partial F}{\partial w}\big(b,u(b),\alpha_2,S''(b)\big)\nabla_{\mathbf{y}}S''(b) = 0,
\end{array}
\end{equation}
where
$$
\begin{array}{rcl}
\nabla_{\mathbf{y}}S''(x) & = & (1-\alpha_1-\alpha_2)\nabla_{\mathbf{y}}W_0''(x)+\alpha_1\nabla_{\mathbf{y}}W_1''(x)+\alpha_2\nabla_{\mathbf{y}}W_2''(x)\medskip\\
& & - (\nabla_{\mathbf{y}}\alpha_1+\nabla_{\mathbf{y}}\alpha_2)W_0''(x)+\nabla_{\mathbf{y}}\alpha_1W_1''(x)+\nabla_{\mathbf{y}}\alpha_2W_2''(x)
\end{array}
$$

System (\ref{CS_Diri_Gradient}) is {\bf linear} in terms of $\nabla_{\mathbf{y}}\alpha_1$ and $\nabla_{\mathbf{y}}\alpha_2$, and it be can easily solved for them.\\

Finally, when using arbitrary BCs as in (\ref{BCs}), $\alpha_1$ and $\alpha_2$ are individualy solutions of the equations in (\ref{BCs_1-2}) in terms of $y_0$ and $y_n$ respectively. Therefore, we have
$$\dfrac{\partial \alpha_1}{\partial y_j} = 0,\forall\,j\neq 0,\quad\textrm{and}\quad\dfrac{\partial \alpha_2}{\partial y_j} = 0,\forall\,j\neq n$$
The remaining partial derivatives can be computed from the following relations obtained by differentiating the equations in (\ref{BCs_1-2}) with respect to $y_0$ and $y_n$ respectively (in other words, using the implicit function theorem)
\begin{equation}\label{BCs_1-2_Gradient}
\begin{array}{l}
\dfrac{\partial g}{\partial u}\big(a,y_0,\alpha_1\big)+\dfrac{\partial g}{\partial v}\big(a,y_0,\alpha_1\big)\dfrac{\partial \alpha_1}{\partial y_0} = 0,\medskip\\
\dfrac{\partial h}{\partial u}(b,y_n,\alpha_2\big)+\dfrac{\partial h}{\partial v}(b,y_n,\alpha_2\big)\dfrac{\partial \alpha_2}{\partial y_n} = 0.
\end{array}
\end{equation}
In summary, let be given $n+1$ knots $\{x_0,\cdots,x_n\}$ and a vector $\mathbf{y} = [y_0,\cdots,y_n]$, then we can compute the cubic and quintic splines $S$ that interpolate the data points $(x_i,y_i)$'s and that correspond to the considered BVP, more precisely, to the BCs of the BVP. Using analytical differentiation, we have the exact spline derivatives, $S',\,S'',\cdots$. Moreover, using automatic differentiation, we have the gradients with respect to the vector $\mathbf{y}$ of $S$ and its derivatives. Therefore, we have all what we need to evaluate the residual of the ODE and its gradient for the obtained spline $S$ at any point. For $\xi\in[a,b]$, the corresponding physics-informed residual is defined by
$$\mathcal{R}(\xi)=F\big(\xi,S(\xi),S'(\xi),S''(\xi)\big).$$
If we denote by $(\xi_j)_j$ the dense set of collocation points, we denote the residual at $\xi_j$ by $\mathcal{R}_j=\mathcal{R}(\xi_j)$, $1\leq j\leq N$, and the residual vector $\mathbf{\mathcal{R}} = \Big(\mathcal{R}_1,\cdots,\mathcal{R}_N\Big)\in\R^N$. The loss function $J$ at the input vector $\mathbf{y}$ is defined to be the $L^2$ or $L^{\infty}$ norm of the residual vector $\mathbf{\mathcal{R}}$. This corresponds exactly to (\ref{ObjectiveFctL2}) for the $L^2$ norm.\\

As example, in the case of Dirichlet BCs, for both cubic and quintic splines, the loss function $J$ can be summarized as follows.
\begin{itemize}
\item Given: the knots $x_0,\cdots,x_n$, the collocation points $\xi_j$'s, and the boundary values $u(a)$ and $u(b)$.
\item Define $J:\mathbf{z}\in\R^{n-1}\longmapsto J(\mathbf{z})\in\R$ by the following
\begin{itemize}
\item    $\mathbf{y} = [u(a), \mathbf{z}, u(b)]$.
\item    Compute the clamped basic spline functions $W_\ell$'s interpolating the points $(x_i,y_i)$ and their derivatives $W_\ell'$, $W_\ell'',\cdots$, and also all of their gradients with respect to $\mathbf{y}$ which depend only on distances $h_i$'s between the knots $x_i$'s.
\item    Solve the nonlinear system for the $\alpha_\ell$'s using Newton-Raphson's algorithm and also compute their gradients.
\item    Compute $\displaystyle S = \sum_{\ell}\alpha_\ell W_\ell$ and all of its derivatives $S'$, $S''$, $\cdots$ along with their gradients.
\item    Evaluate, for all $1\leq j\leq N$, the physics-informed residual $\mathcal{R}_j$ and its gradient at the collocation point $\xi_j$
    	$$\mathcal{R}_j = F\Big(\xi_j,S(\xi_j),S'(\xi_j),S''(\xi_j)\Big).$$
\item    Return the $L^2$ or $L^{\infty}$ norm of the physics-informed residual vector $\mathbf{\mathcal{R}}$ as $J(\mathbf{y})$ and its gradient $\nabla_{\mathbf{y}} J(\mathbf{y})$.
\end{itemize}
\end{itemize}

SPINS consists on minimizing the loss function $J$ in order to provide the optimal vector $\mathbf{y}$ from which we will generate the optimal spline that verifies the BCs of the BVP and, at the same time, makes the physics-informed residual minimal.\\

To optimize the spline coefficients, we utilized the L-BFGS-B algorithm, leveraging the exact analytical gradients provided by SPINS. Since the L-BFGS-B algorithm is iterative, so it should be initialized. However, we do not have any information about the unknown vector $\mathbf{y}$. We decided to initialize with arbitrary random vectors. An alternative to this initialization could be by using traditional numerical methods but with coarse subdivisions of $[a,b]$ and accept rough estimates of the derivatives.

\section{Numerical Simulation}
In order to validate SPINS, we present in this section several BVPs with a wide variability in the nonlinearity of both the ODE and the BCs in addition to the choice of various exact solutions.\\

The knots are chosen to be either uniformly distributed between $a$ and $b$, i.e. $x_i = (b-a)\dfrac{i}{n}$, $i=0,\cdots,n$, or the Chebeychev points. The initial vector $\mathbf{y} = [y_0,\cdots,y_n]$ is considered randomly between $-1$ and $1$. When dealing with Dirichlet BCs we take $y_0= u(a)$ and $y_n=u(b)$.\\

The function $u$ denotes the exact solution of the BVP, and $S_3$ and $S_5$ respectively the cubic and quintic splines solutions obtained by SPINS. Moreover, we denote by $S_{3,exact}$ and $S_{5,exact}$ respectively the cubic and quintic splines obtained by interpolating the exact data points $(x_i,u(x_i))$, $i=0,\cdots,n$ that lie on the curve of $u$, but they satisfy the BCs of the BVP.\\

The physics-informed residual is evaluated at $100$ collocation points uniformly distributed between $a$ and $b$. Its $L^2$ norm is computed according to Equation (\ref{ObjectiveFctL2}) (and applying a square root). The relative $L^2$, $H^1$ and $H^2$ errors are computed as a percentage according to the following formula
$$\textrm{Error}_{L^2} = 100\dfrac{\|u-S\|_{L^2}}{\|u\|_{L^2}},\quad
\textrm{Error}_{H^1} = 100\dfrac{\|u-S\|_{H^1}}{\|u\|_{H^1}},\quad
\textrm{Error}_{H^2} = 100\dfrac{\|u-S\|_{H^2}}{\|u\|_{H^2}},$$
that are also evaluated at the collocation points, where $\|u\|_{H^1}^2 = \|u\|_{L^2}^2+\|u'\|_{L^2}^2$ and $\|u\|_{H^2}^2 = \|u\|_{L^2}^2+\|u'\|_{L^2}^2+\|u''\|_{L^2}^2$.

\subsection{Example 1: Cubic SPINS and Dirichlet BCs}
Consider the following BVP with Dirichlet BCs
\begin{equation}\label{BVP1}
\left\{\begin{array}{l}
u''(x)+u'(x)+\sin\big(u'(x)u(x)\big) = f(x),\quad 0\leq x\leq 2\pi,\medskip\\
u(0) = u(2\pi) = 0,
\end{array}\right.
\end{equation}
with exact solution $u(x) =\sin(x)$. The right hand side $f$ is obtained by plugging the exact solution into the left hand side of the ODE.
We solved BVP (\ref{BVP1}) using cubic SPINS where the knots are considered uniformly distributed.
\begin{table}[H]
    \centering
    \caption{Cubic SPINS for BVP (\ref{BVP1}) (Dirichlet BCs).}
    \label{tab:Example1}
{\scriptsize{
\begin{tabular}{|c||c|c|c|c||c|c|c|c|}
\toprule[1.5pt]
\multicolumn{1}{c}{} & \multicolumn{4}{c}{$S_{3,exact}$} & \multicolumn{4}{c}{$S_3$ (SPINS)}\\
 \cmidrule[1.5pt](lr){2-5}\cmidrule[1.5pt](lr){6-9}
\multicolumn{1}{|c||}{Knots} & 
\multicolumn{1}{c}{Residual} & \multicolumn{1}{c}{Error$_{L^2}$} & 
\multicolumn{1}{c}{Error$_{H^1}$} & \multicolumn{1}{c||}{Error$_{H^2}$} & 
\multicolumn{1}{c}{Residual} & \multicolumn{1}{c}{Error$_{L^2}$} & 
\multicolumn{1}{c}{Error$_{H^1}$} & \multicolumn{1}{c|}{Error$_{H^2}$}\\
\hline
$5$ & 9.29e-02 & 2.07e+00 \% & 3.48e+00 \% & 8.04e+00 \% & 9.21e-02 & 2.85e+00 \% & 3.98e+00 \% & 8.57e+00 \%\\
$8$ & 2.36e-02 & 1.36e-01 \% & 3.81e-01 \% & 1.93e+00 \% & 2.36e-02 & 1.70e-01 \% & 3.87e-01 \% & 1.92e+00 \%\\
$10$ & 1.42e-02 & 4.55e-02 \% &1.62e-01 \% & 1.16e+00 \% & 1.42e-02 & 8.62e-02 \% & 1.73e-01 \% & 1.16e+00 \%\\
\hline
\bottomrule[1.5pt]
\end{tabular}}}
\end{table}

\begin{figure}[!htbp]
    \centering
    \includegraphics[width=0.8\textwidth]{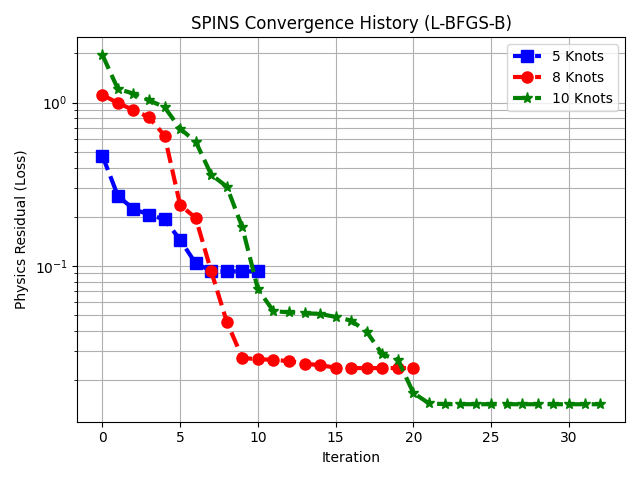}
    \caption{Convergence history of the cubic SPINS for BVP (\ref{BVP1}). The physics-informed residual is plotted on a logarithmic scale, demonstrating stable descent to high-precision levels using L-BFGS-B.}
    \label{fig:example1_Loss}
\end{figure}

\begin{figure}[!htbp]
    \centering
    \includegraphics[width=0.48\textwidth]{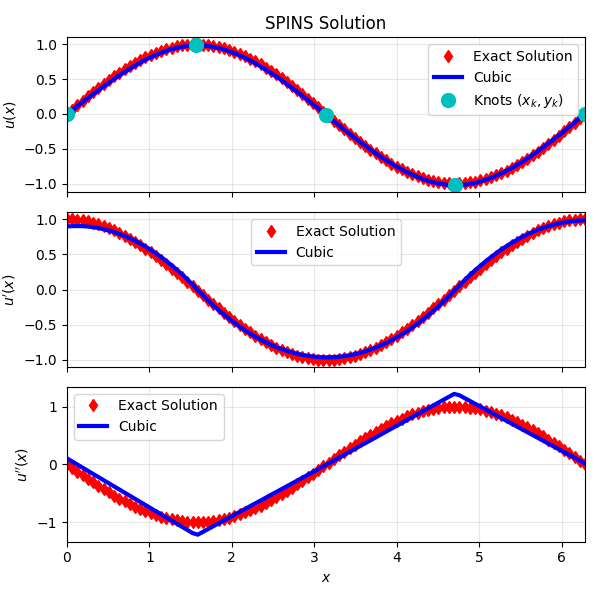}
    \includegraphics[width=0.48\textwidth]{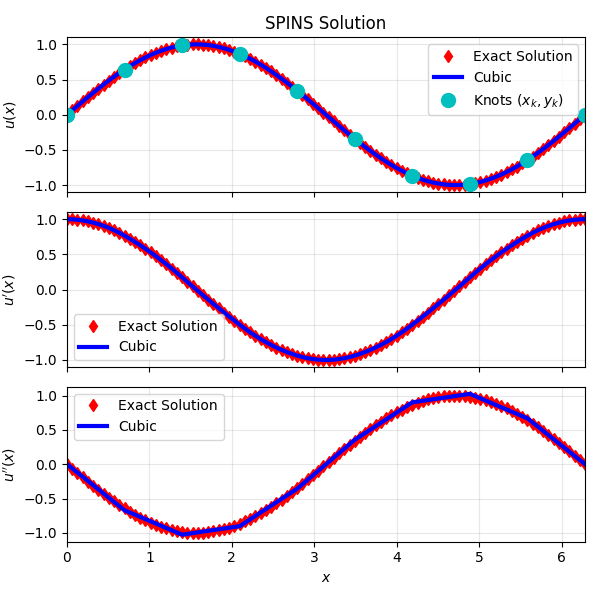}
    \caption{The cubic SPINS and the exact solutions with their first and second order derivatives for $5$ knots (left) 
    and $10$ knots (right) for BVP (\ref{BVP1}).}
    \label{fig:example1_Solution}
\end{figure}

\subsection{Example 2: Quintic SPINS and Dirichlet BCs}
Consider the following BVP with Dirichlet BCs
\begin{equation}\label{BVP2}
\left\{\begin{array}{l}
u''(x)+u'(x)^2-u'(x)^3+u'(x)u(x)+\sin(u(x))^3 = f(x),\quad 0\leq x\leq \pi,\medskip\\
u(0) = u(\pi) = 0,
\end{array}\right.
\end{equation}
with exact solution $u(x) =xe^{-x}\sin x$. 
We solved BVP (\ref{BVP2}) using quintic SPINS where the knots are considered uniformly distributed.

\begin{table}[H]
    \centering
    \caption{Quintic SPINS for BVP (\ref{BVP2}) (Dirichlet BCs).}
    \label{tab:Example2}
{\scriptsize{
\begin{tabular}{|c||c|c|c|c||c|c|c|c|}
\toprule[1.5pt]
\multicolumn{1}{c}{} & \multicolumn{4}{c}{$S_{5,exact}$} & \multicolumn{4}{c}{$S_5$ (SPINS)}\\
 \cmidrule[1.5pt](lr){2-5}\cmidrule[1.5pt](lr){6-9}
\multicolumn{1}{|c||}{Knots} & 
\multicolumn{1}{c}{Residual} & \multicolumn{1}{c}{Error$_{L^2}$} & 
\multicolumn{1}{c}{Error$_{H^1}$} & \multicolumn{1}{c||}{Error$_{H^2}$} & 
\multicolumn{1}{c}{Residual} & \multicolumn{1}{c}{Error$_{L^2}$} & 
\multicolumn{1}{c}{Error$_{H^1}$} & \multicolumn{1}{c|}{Error$_{H^2}$}\\
\hline
$5$ & 5.59e-03 & 1.51e-01 \% & 4.30e-01 \% & 9.42e-01 \% & 5.41e-03 & 1.09e+00 \% & 1.10e+00 \% & 1.09e+00 \%\\
$8$ & 4.89e-04 & 4.03e-03 \% & 2.01e-02 \% & 7.95e-02 \% & 4.87e-04 & 1.76e-02 \% & 2.58e-02 \% & 7.96e-02 \%\\
$10$ & 1.60e-04 & 7.67e-04 \%  & 4.97e-03 \% & 2.57e-02 \% & 1.59e-04 & 1.23e-03 \% & 5.22e-03 \% & 2.57e-02 \%\\
\hline
\bottomrule[1.5pt]
\end{tabular}}}
\end{table}

\begin{figure}[!htbp]
    \centering
    \includegraphics[width=0.48\textwidth]{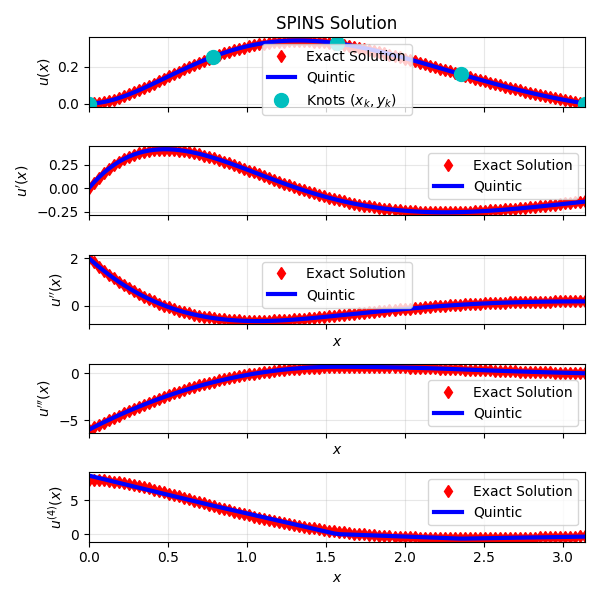}
    \includegraphics[width=0.48\textwidth]{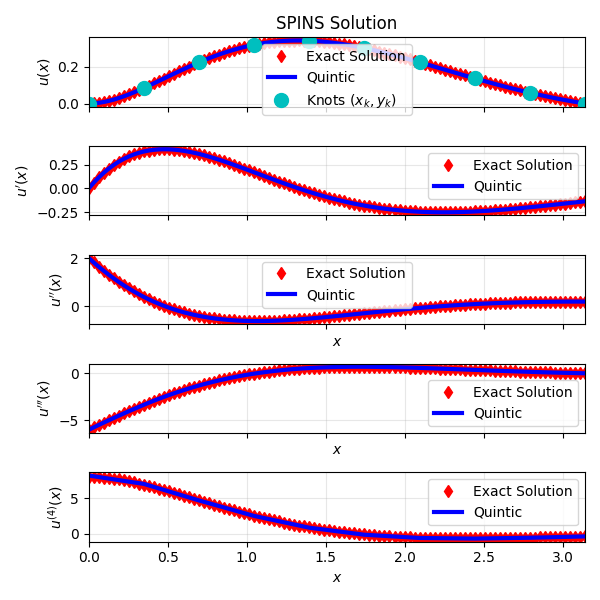}
    \caption{The quintic SPINS and the exact solutions with their derivatives up to order four for $5$ knots (left)
    and $10$ knots (right) for BVP (\ref{BVP2}).}
    \label{fig:example2_Solution}
\end{figure}
  
\subsection{Example 3: Quintic SPINS and Neumann BCs}
Consider the following BVP with Neumann BCs
\begin{equation}\label{BVP3}
\left\{\begin{array}{l}
-u''(x)+\dfrac{1}{1+u''(x)^2}+u(x)^3 = f(x),\quad 0\leq x\leq 2,\medskip\\
u'(0) = \pi,\,\,u'(2) = \dfrac{\pi}{3},
\end{array}\right.
\end{equation}
with exact solution $u(x) =\dfrac{\sin(\pi x)}{1+x}$.
We solved BVP (\ref{BVP3}) using quintic SPINS where the knots are considered to be the Chebeychev points on the interval $[0,2]$.
\begin{table}[H]
    \centering
    \caption{Quintic SPINS for BVP (\ref{BVP3}) (Neumann BCs + Chebeychev knots).}
    \label{tab:Example3}
{\scriptsize{
\begin{tabular}{|c||c|c|c|c||c|c|c|c|}
\toprule[1.5pt]
\multicolumn{1}{c}{} & \multicolumn{4}{c}{$S_{5,exact}$} & \multicolumn{4}{c}{$S_5$ (SPINS)}\\
 \cmidrule[1.5pt](lr){2-5}\cmidrule[1.5pt](lr){6-9}
\multicolumn{1}{|c||}{Knots} & 
\multicolumn{1}{c}{Residual} & \multicolumn{1}{c}{Error$_{L^2}$} & 
\multicolumn{1}{c}{Error$_{H^1}$} & \multicolumn{1}{c||}{Error$_{H^2}$} & 
\multicolumn{1}{c}{Residual} & \multicolumn{1}{c}{Error$_{L^2}$} & 
\multicolumn{1}{c}{Error$_{H^1}$} & \multicolumn{1}{c|}{Error$_{H^2}$}\\
\hline
$5$ & 8.63e-02 & 7.46e-01 \% & 1.05e+00 \% & 1.80e+00 \% & 7.25e-02 & 6.05e+00 \% & 1.89e+00 \% & 1.59e+00 \%\\ 
$8$ & 1.34e-02 & 4.03e-02 \% & 9.11e-02 \% & 2.53e-01 \% & 1.22e-02 & 2.70e-01 \% & 1.15e-01 \% & 2.48e-01 \%\\ 
$10$ & 3.72e-03 & 5.83e-03 \% & 1.74e-02 \% & 7.07e-02 \% & 3.46e-03 & 1.94e-01 \% & 6.19e-02 \% & 7.65e-02 \%\\ 
\hline
\bottomrule[1.5pt]
\end{tabular}}}
\end{table}

\begin{figure}[!htbp]
    \centering
    \includegraphics[width=0.48\textwidth]{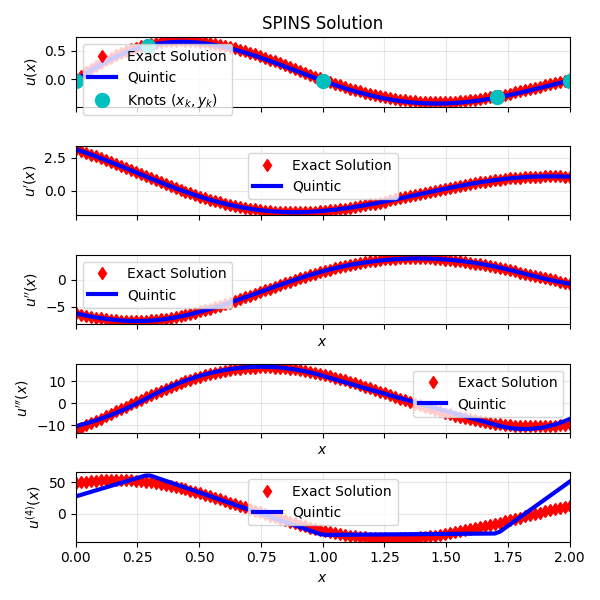}
    \includegraphics[width=0.48\textwidth]{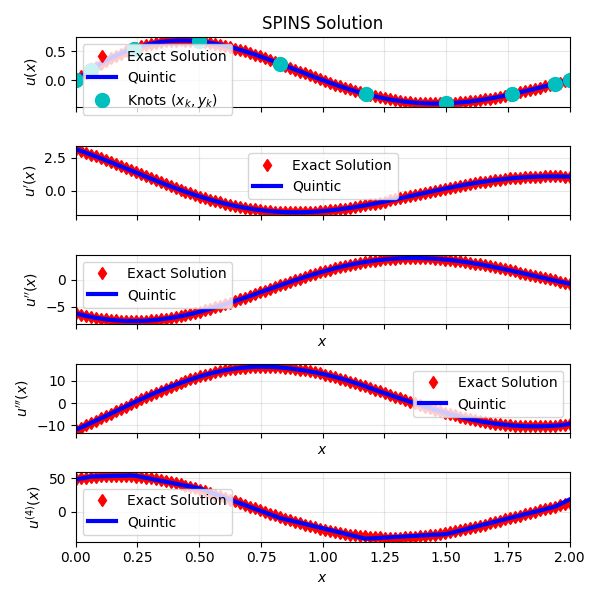}
    \caption{The quintic SPINS and the exact solutions with their derivatives up to order four for $5$ knots (left)
    and $10$ knots (right) for BVP (\ref{BVP3}).}
    \label{fig:example3_Solution}
\end{figure}
\subsection{Example 4: Cubic and Quintic SPINS with Robin BCs}
Consider the following BVP with Robin BCs
\begin{equation}\label{BVP4}
\left\{\begin{array}{l}
u''(x)+u'(x)^2+\sin\big(u'(x)u(x)\big) = f(x),\quad -1\leq x\leq 1,\medskip\\
u(-1)+u'(-1) = 1,\,\,u(1)-u'(1) = 1,
\end{array}\right.
\end{equation}
with exact solution $u(x) =\dfrac{1}{1+x^2}$.
We solved BVP (\ref{BVP4}) using both cubic and quintic SPINS where the knots are considered uniformly distributed on the interval $[-1,1]$.
\begin{table}[H]
    \centering
    \caption{Cubic SPINS for BVP (\ref{BVP4}) (Robin BCs).}
    \label{tab:Example4}
{\scriptsize{
\begin{tabular}{|c||c|c|c|c||c|c|c|c|}
\toprule[1.5pt]
\multicolumn{1}{c}{} & \multicolumn{4}{c}{$S_{3,exact}$} & \multicolumn{4}{c}{$S_3$ (SPINS)}\\
 \cmidrule[1.5pt](lr){2-5}\cmidrule[1.5pt](lr){6-9}
\multicolumn{1}{|c||}{Knots} & 
\multicolumn{1}{c}{Residual} & \multicolumn{1}{c}{Error$_{L^2}$} & 
\multicolumn{1}{c}{Error$_{H^1}$} & \multicolumn{1}{c||}{Error$_{H^2}$} & 
\multicolumn{1}{c}{Residual} & \multicolumn{1}{c}{Error$_{L^2}$} & 
\multicolumn{1}{c}{Error$_{H^1}$} & \multicolumn{1}{c|}{Error$_{H^2}$}\\
\hline
$5$ & 9.79e-02 & 1.83e-01 \% & 1.03e+00 \% & 6.96e+00 \% & 9.76e-02 & 1.44e+00 \% & 1.87e+00 \% & 7.10e+00 \%\\ 
$8$ & 4.19e-02 & 2.68e-02 \% & 2.61e-01 \% & 2.98e+00 \% & 4.19e-02 & 1.29e-01 \% & 2.89e-01 \% & 2.99e+00 \%\\ 
$10$ & 2.40e-02 & 7.85e-03 \% & 9.99e-02 \% & 1.71e+00 \% & 2.40e-02 & 4.50e-02 \% & 1.10e-01 \% & 1.71e+00 \%\\ 
\hline
\bottomrule[1.5pt]
\end{tabular}}}
\end{table}

\begin{figure}[!htbp]
    \centering
    \includegraphics[width=0.48\textwidth]{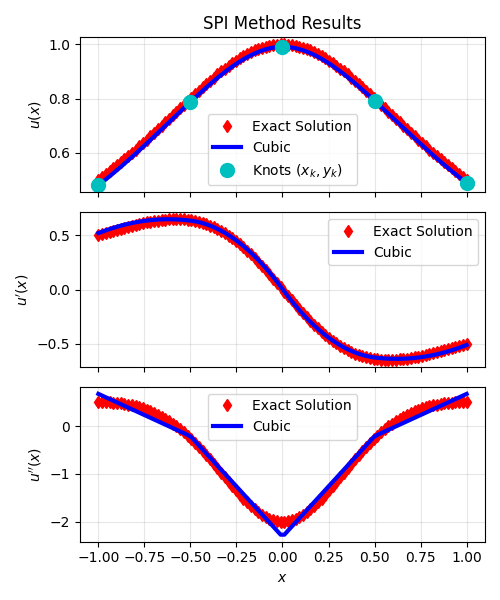}
    \includegraphics[width=0.48\textwidth]{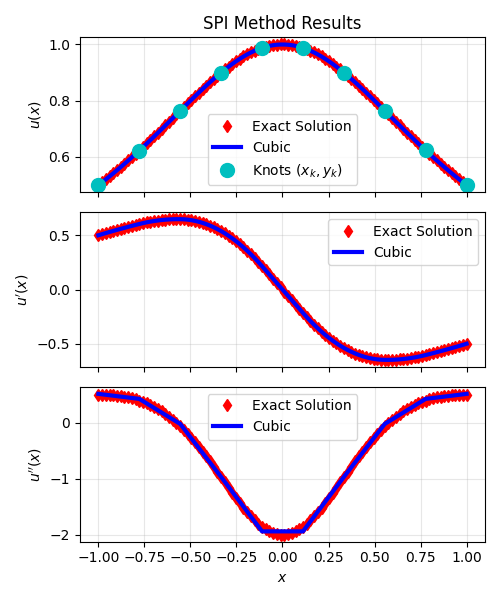}
    \caption{The cubic SPINS and the exact solutions with their first and second order derivatives for $5$ knots (left)
    and $10$ knots (right) for BVP (\ref{BVP4}).}
    \label{fig:example4_Solution}
\end{figure}

\begin{table}[H]
    \centering
    \caption{Quintic SPINS for BVP (\ref{BVP4}) (Robin BCs).}
    \label{tab:Example4bb}
{\scriptsize{
\begin{tabular}{|c||c|c|c|c||c|c|c|c|}
\toprule[1.5pt]
\multicolumn{1}{c}{} & \multicolumn{4}{c}{$S_{5,exact}$} & \multicolumn{4}{c}{$S_5$ (SPINS)}\\
 \cmidrule[1.5pt](lr){2-5}\cmidrule[1.5pt](lr){6-9}
\multicolumn{1}{|c||}{Knots} & 
\multicolumn{1}{c}{Residual} & \multicolumn{1}{c}{Error$_{L^2}$} & 
\multicolumn{1}{c}{Error$_{H^1}$} & \multicolumn{1}{c||}{Error$_{H^2}$} & 
\multicolumn{1}{c}{Residual} & \multicolumn{1}{c}{Error$_{L^2}$} & 
\multicolumn{1}{c}{Error$_{H^1}$} & \multicolumn{1}{c|}{Error$_{H^2}$}\\
\hline
$5$ & 4.32e-03 & 5.59e-03 \% & 3.85e-02 \% & 3.07e-01 \% & 4.27e-03 & 1.81e-02 \% & 3.90e-02 \% & 3.04e-01 \%\\ 
$8$ & 7.61e-03 & 6.74e-03 \% & 6.35e-02 \% & 5.39e-01 \% & 7.43e-03 & 1.86e-02 \% & 6.05e-02 \% & 5.28e-01 \%\\ 
$10$ & 1.62e-03 & 7.27e-04 \% & 9.12e-03 \% & 1.16e-01 \% & 1.61e-03 & 8.04e-04 \% & 8.53e-03 \% & 1.15e-01 \%\\
\hline
\bottomrule[1.5pt]
\end{tabular}}}
\end{table}

\begin{figure}[!htbp]
    \centering
    \includegraphics[width=0.48\textwidth]{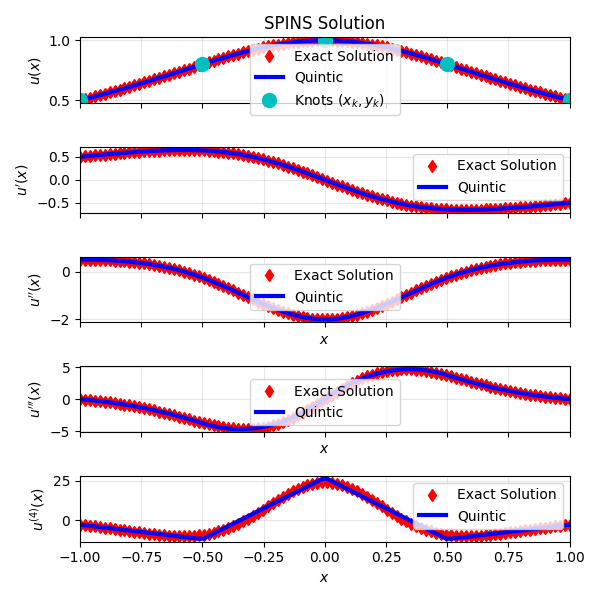}
    \includegraphics[width=0.48\textwidth]{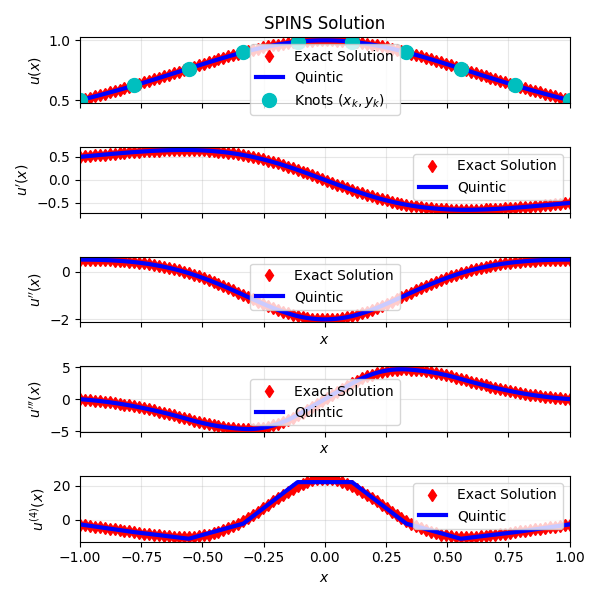}
    \caption{The quintic SPINS and the exact solutions with their derivatives up to order four for $5$ knots (left)
    and $10$ knots (right) for BVP (\ref{BVP4}).}
    \label{fig:example4bb_Solution}
\end{figure}
\subsection{Example 5: Quintic SPINS and Nonlinear BCs}
Consider the following BVP with nonlinear BCs
\begin{equation}\label{BVP5}
\left\{\begin{array}{l}
-u''(x)+u'(x)+\sin\big(u(x)\big) = f(x),\quad -\pi\leq x\leq \pi,\medskip\\
\sqrt{1+u(-\pi)^2}+u'(-\pi) = 1+e^{\pi}\pi\approx 73.6986,\,\,u(\pi)+u'(\pi)^3 =  -e^{-3\pi}\pi^3 \approx -0.0025,
\end{array}\right.
\end{equation}
with exact solution $u(x) =xe^{-x}\sin x$. We solved BVP (\ref{BVP5}) using quintic SPINS where the knots are considered uniformly distributed on the interval $[-\pi,\pi]$.
\begin{table}[H]
    \centering
    \caption{Quintic SPINS for BVP (\ref{BVP5}) (Nonlinear BCs).}
    \label{tab:Example5}
{\scriptsize{
\begin{tabular}{|c||c|c|c|c||c|c|c|c|}
\toprule[1.5pt]
\multicolumn{1}{c}{} & \multicolumn{4}{c}{$S_{5,exact}$} & \multicolumn{4}{c}{$S_5$ (SPINS)}\\
 \cmidrule[1.5pt](lr){2-5}\cmidrule[1.5pt](lr){6-9}
\multicolumn{1}{|c||}{Knots} & 
\multicolumn{1}{c}{Residual} & \multicolumn{1}{c}{Error$_{L^2}$} & 
\multicolumn{1}{c}{Error$_{H^1}$} & \multicolumn{1}{c||}{Error$_{H^2}$} & 
\multicolumn{1}{c}{Residual} & \multicolumn{1}{c}{Error$_{L^2}$} & 
\multicolumn{1}{c}{Error$_{H^1}$} & \multicolumn{1}{c|}{Error$_{H^2}$}\\
\hline
$5$ & 1.16e+00 & 2.87e+00 \% & 3.08e+00 \% & 2.52e+00 \% & 9.50e-01 & 5.11e+00 \% & 3.12e+00 \% & 2.20e+00 \%\\ 
$8$ & 1.37e-01 & 1.03e-01 \% & 1.89e-01 \% & 3.04e-01 \% & 1.11e-01 & 9.23e-01 \% & 4.72e-01 \% & 2.96e-01 \%\\ 
$10$ & 4.24e-02 & 1.68e-02 \% & 4.10e-02 \% & 9.50e-02 \% & 3.82e-02 & 9.18e-02 \% & 6.51e-02 \% & 8.67e-02 \%\\ 
\hline
\bottomrule[1.5pt]
\end{tabular}}}
\end{table}

\begin{figure}[!htbp]
    \centering
    \includegraphics[width=0.48\textwidth]{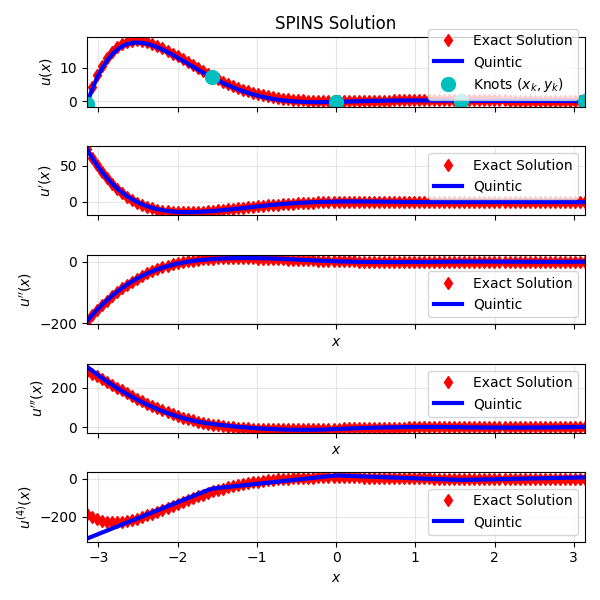}
    \includegraphics[width=0.48\textwidth]{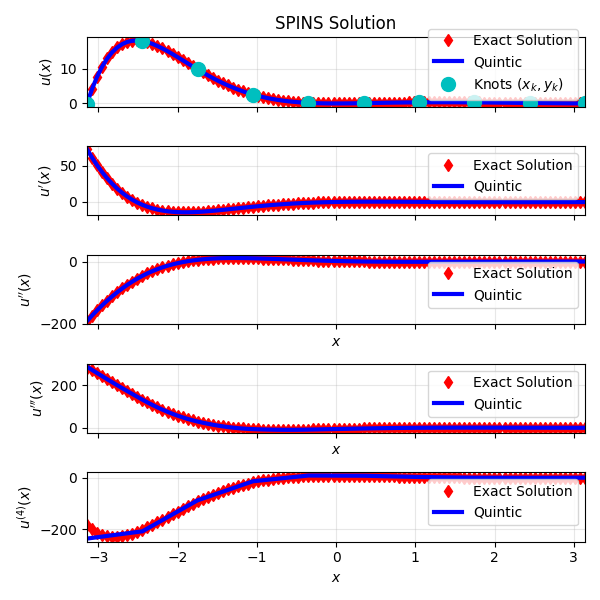}
    \caption{The quintic SPINS and the exact solutions with their derivatives up to order four for $5$ knots (left)
    and $10$ knots (right) for BVP (\ref{BVP5}).}
    \label{fig:example5_Solution}
\end{figure}

\begin{figure}[!htbp]
    \centering
    \includegraphics[width=0.43\textwidth]{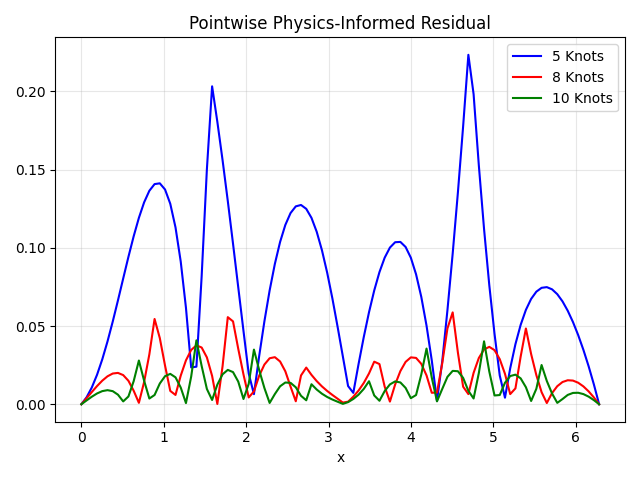}
	\includegraphics[width=0.43\textwidth]{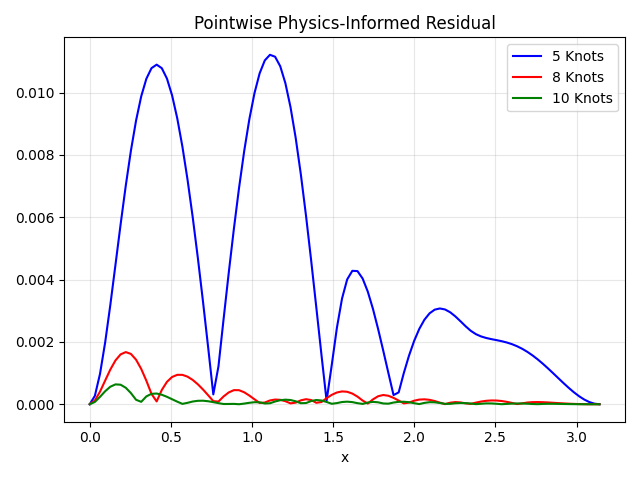}
    \includegraphics[width=0.43\textwidth]{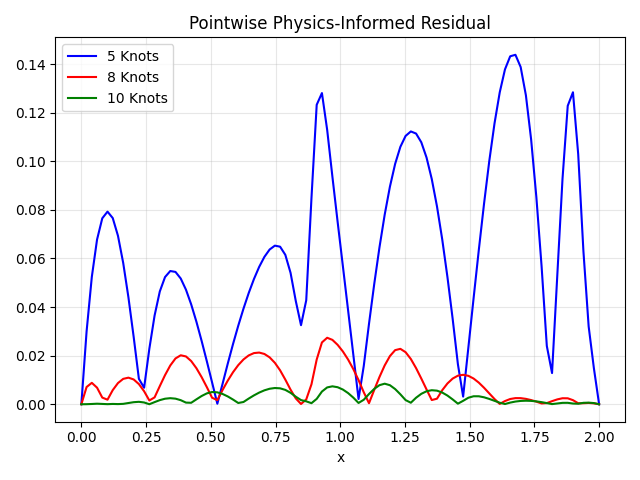}
    \includegraphics[width=0.43\textwidth]{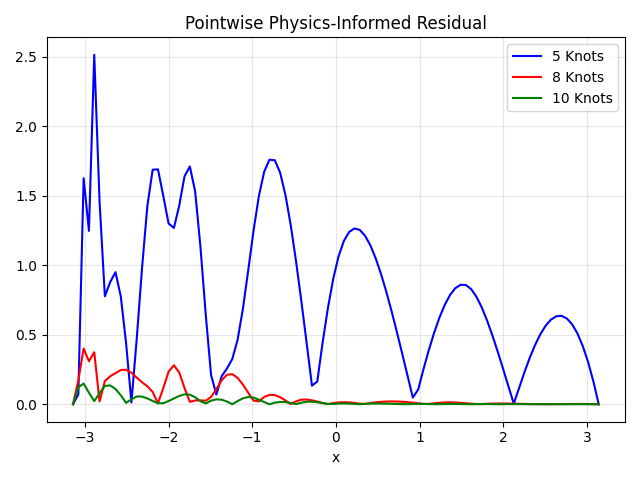}
    \caption{The pointwise physics-informed residual on the interval $[a,b]$ at the last iteration of the minimization process (after convergence) for BVPs (\ref{BVP1}), (\ref{BVP2}), (\ref{BVP3}) and (\ref{BVP5}).}
    \label{fig:examples_ODE_residual}
\end{figure}
\section{Discussion of the Numerical Results}\label{Interpretation}
SPINS provides a spline solution to a BVP by minimizing the $L^2$ norm of the physics-informed residual evaluated at the collocation points. At each iteration of the minimization process, the spline is obtained by interpolating the data points $(x_i,y_i)$. In general, to compute cubic or quintic splines, additional information on the splines derivatives at the extreme knots are required. In our case, if they are not provided as explicit boundary conditions, these derivatives values are extracted from the ODE itself, and in the fully nonlinear BVPs, this requires solving nonlinear equations of one real variable, or systems of at most $4$ equations (in the worst scenario).\\

Since we are searching for splines that will approximate the exact solution $u$, then in addition to the SPINS solution $S$, we computed for validation the splines $S_{exact}$ that we obtain by interpolating the exact data points $(x_i,u(x_i))_i$'s with the BCs as we do for $S$. Also we evaluated their corresponding physics-informed residual and the $L^2$, $H^1$ and $H^2$ norms of the relative error between $u$ and $S_{exact}$ and between $u$ and $S$.\\

As expected, the relative errors of $S_{exact}$ are in general smaller than the relative errors of $S$, and this is because $S_{exact}$ interpolates the values of the exact solution $u$. However, in almost all of the numerical simulations, we noticed that the $L^2$ norm of the physics-informed residual of the SPINS solution $S$ is smaller than the $L^2$ norm of the physics-informed residual of the exact approximating spline $S_{exact}$, see Tables \ref{tab:Example1}-\ref{tab:Example5}. This means that, although $S_{exact}$ interpolates the values of the exact solution $u$, $S_{exact}$ is not the optimal solution for the minimization of the physics-informed residual. And SPINS provides the optimal spline that minimizes the physics-informed residual.\\

On the other hand, the relative errors of the SPINS solution decrease when the number of interpolating knots increases and also when the order $q$ of the spline increases (from cubic to quintic) which is expected in numerical analysis. But the additional remark here is that our observations are not restricted to the $L^2$ norm of the relative errors but they are also valid for the $H^1$ and $H^2$ norms (and also for the $H^3$ and $H^4$ norms for the case of the quintic SPINS). In other words, using cubic or quintic splines to approximate the exact solution  of a BVP by simply minimizing the $L^2$ norm of the corresponding physics-informed residual, provides a smooth solution whose derivatives also approximate the derivatives of the exact solution, see Figures \ref{fig:example1_Solution}, \ref{fig:example2_Solution}, \ref{fig:example3_Solution}, \ref{fig:example4_Solution}, \ref{fig:example4bb_Solution}, and \ref{fig:example5_Solution}.\\

As for the pointwise physics-informed residual, it is illustrated on Figure \ref{fig:examples_ODE_residual} on the interval $[a,b]$ at the last iteration of the minimization process (after convergence of the L-BFGS-B) for BVPs (\ref{BVP1}), (\ref{BVP2}), (\ref{BVP3}) and (\ref{BVP5}).
We see on Figure \ref{fig:examples_ODE_residual} the space variation of this residual and we can locate the regions of high residual. This may help in increasing the number of knots near these regions to compensate the error. Furthermore, we notice that for $n=10$ knots, the overall spatial residual distribution is very small compared to the case with $n=8$ knots which is very small compared to the case with $n=5$ knots.\\

Figure \ref{fig:example1_Loss} shows, for BVP (\ref{BVP1}), the curves of the $L^2$ norm of the physics-informed residual (the loss history) versus the iterations numbers of the L-BFGS-B algorithm. First, it is important to notice that the convergence of the gradient-based optimizer L-BFGS-B used for SPINS takes very few iterations: about $10$ iterations for $5$ knots, $20$ iterations for $8$ knots and $32$ iterations for $10$ knots. This number of iterations is negligible compared to the hundreds or thousands of iterations (epochs) required using Adam optimizer for PINNS. Most probably, this is due to the fact that with SPINS we are minimizing an objective function defined on very small dimensional spaces (e.g., $\R^5$ , $\R^8$ or $\R^{10}$) while in PINNS, the vector of unknown parameters (the weights and the biases of the neural network) belong to very high dimensional (hundreds or more) spaces depending on the number of hidden layers and of neurons in each layer. Second, we noticed that for large number of knots, the optimization takes much more iterations to converge than the case of small number of knots. And here we mean by convergence that the optimization algorithm stabilizes and no more improvements to the solution can be provided. Actually, this happens when the loss value becomes close to (and as mentioned above, in most of the examples, less than) the loss value of the exact spline $S_{exact}$. Indeed, this is due to the fact that choosing a random initial state $\mathbf{y}$ with a large number of knots introduces lots of artificial variations in the spline solution and the optimization process requires more time to take rid of these unrealistic variations. However, the loss value for larger number of knots is much smaller than the loss value of smaller number of knots. This gives us an idea about the rate of convergence of SPINS.\\

The theoretical foundation for spline-based approximation is characterized by optimal error bounds in both the $L^2$ and $L^\infty$ norms. For a target function $u \in H^{q+1}[a, b]$, where $q$ denotes the spline degree, the approximation $S_{exact}$ and its $r$-th derivatives satisfy the general bound $\|u^{(r)} - S_{exact}^{(r)}\|_{L^p} \leq C_{r,p} h^{q+1-r} \|u^{(q+1)}\|_{L^p}$ for $p \in \{2, \infty\}$. Specifically, for cubic splines ($q=3$), the function approximation converges at an optimal rate of $O(h^4)$ in both norms, while the first and second-order derivatives converge at $O(h^3)$ and $O(h^2)$ respectively, \cite{douglas1975, hall1976}. In contrast, the quintic splines ($q=5$) significantly enhance these rates, yielding $O(h^6)$ convergence for the solution, $O(h^5)$ for the first-order differential operator and $O(h^4)$ for the second-order differential operator, \cite{schultz1973}.\\

We estimated the rate of convergence, denoted by $s$,  provided by SPINS for BVP (\ref{BVP1}) which uses cubic splines and BVP (\ref{BVP2}) which uses quintic splines according to the following formula
$$\|u-S\| \approx C h^s \|u\|$$
The value of $s$ is estimated as follows
$$ s  \approx \dfrac{\ln(Error_{10})-\ln(Error_{5})}{\ln(h_{10})-\ln(h_{5})},$$
where $h_n = \dfrac{b-a}{n-1}$, $Error_n$ is each of the $L^2$ norm of the physics-informed residual, $L^2$, $H^1$ and $H^2$ norms of the errors, and $n$ is the number of knots.
\begin{table}[H]
    \centering
    \caption{Rate of convergence $s$ for BVPs (\ref{BVP1}) and (\ref{BVP2}).}
    \label{tab:OrderAccuracy}
{\scriptsize{
\begin{tabular}{|c|c|c|c|c|}
\toprule[1.5pt]
\multicolumn{1}{c}{} & \multicolumn{2}{c}{Cubic} & \multicolumn{2}{c}{Quintic}\\
 \cmidrule[1.5pt](lr){2-3}\cmidrule[1.5pt](lr){4-5}
\multicolumn{1}{c}{} & 
\multicolumn{1}{c}{$S_{3,exact}$} & \multicolumn{1}{c}{$S_3$ (SPINS)} & 
\multicolumn{1}{c}{$S_{5,exact}$} & \multicolumn{1}{c}{$S_5$ (SPINS)}\\
\hline
Residual$_{L^2}$ & 2.32 & 2.31 & 4.38 & 4.35\\ 
\hline
\hline
Error$_{L^2}$ & 4.71 & 4.31 & 6.51 & 8.37\\ 
Error$_{H^1}$ & 3.78 & 3.87 & 5.50 & 6.60\\ 
Error$_{H^2}$ & 2.39 & 2.47 & 4.44 & 4.62\\ 
\hline
\bottomrule[1.5pt]
\end{tabular}}}
\end{table}

First of all, for the case of $S_{exact}$, the numerical estimates of $s$ shown in Table \ref{tab:OrderAccuracy} validate the theoretical results \cite{schultz1973, douglas1975, hall1976} for both the cubic and quintic splines.\\
Second, regarding the $L^2$ norm of the physics-informed residual, we have almost the same rate of convergence for $S_{exact}$ and the SPINS solution $S$ in each of the cubic and quintic cases, where it is about $2$ for the cubic case and $4$ for the quintic case.\\
Moreover, Table \ref{tab:OrderAccuracy} shows that, the rate of convergence of the {\bf cubic} spline provided by SPINS is very close to the rate of convergence obtained from the exact cubic spline for all measured quantities. It is greater than $4$ for the $L^2$ norm of the error, greater than $3$ for the $H^1$ norm of the error, and greater than $2$ for the $H^2$ norm of the error.\\
However, while the estimates of the exact quintic spline $S_{5,exact}$ validate the theoretical results, we see that 
the rate of convergence of the quintic spline $S_5$ provided by SPINS is much better. Indeed, it is almost the {\bf double} of the rate of convergence of the cubic SPINS. It is greater than $8$ for the $L^2$ norm of the error, greater than $6$ for the $H^1$ norm of the error, and greater than $4$ for the $H^2$ norm of the error. And this shows the efficacy of quintic SPINS in all cases, where it provides highly accurate approximations.\\

\begin{figure}[!htbp]
    \centering
    \includegraphics[width=0.48\textwidth]{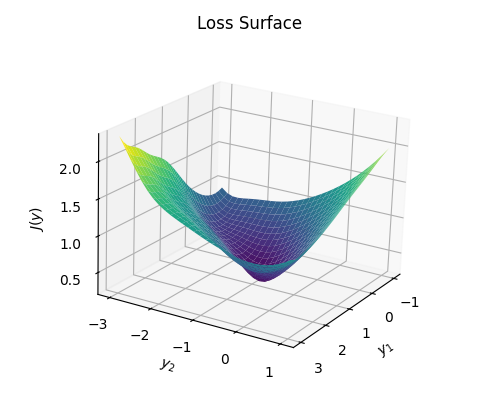}
	\includegraphics[width=0.48\textwidth]{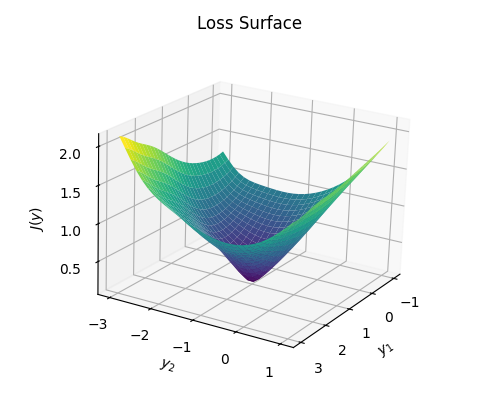}
    \caption{The loss surfaces for BVP (\ref{BVP1}) using cubic SPINS with $4$ knots. Left: the knots are uniformly distributed. Right: The knots are at the local extrema of the exact solution $u$.}
    \label{fig:LossSurfaces}
\end{figure}

To have an idea about the shape of the loss function $J$, we considered BVP (\ref{BVP1}) with Dirichlet BCs and applied cubic SPINS with $4$ knots only, $\{x_0 = 0,x_1,x_2,x_3=2\pi\}$. So $J$ is function of only two variables $\mathbf{y} = [y_1,y_2]$ since $y_0 = u(0)$ and $y_3 = u(2\pi)$. Figure \ref{fig:LossSurfaces} shows the surfaces of $J$ for two different choices of the interior knots: uniformly distributed (left) and at the local extrema (right). For this case study, and unlike the loss function in the case of PINNs, $J$ has a unique global minimum on the considered domain $[-1,3]\times[-3,1]$. It is clear on this figure that choosing the knots at the local extrema of $u$ gives a better shape for the surface of $J$ where the optimization process is much faster. This observation for this particular simple scenario explains the efficiency of the L-BFGS-B optimizer when used for SPINS in all of the above simulations.

\section{Strategies for Better Results}\label{Strategies}
\subsection{Strategy I: From Small to Larger Numbers of Knots}
According to the previous numerical results, the accuracy of the numerical solution with larger number of knots (for instance $10$ knots) is in general much better than the accuracy with small number of knots (for instance $5$ knots), which is natural in numerical analysis. However, as mentioned in Section \ref{Interpretation}, with large number of knots, the optimization takes longer to converge than the case of small number of knots. On the other hand, when we start the optimization process with a large number of knots but initialized with a non-random state $\mathbf{y}$ that is close to the unkown exact solution, the convergence is very fast.\\
To benefit from the fast convergence when choosing a random initial state with a small number of knots, and from the fast convergence when using an initial state close to the exact solution with a large number of knots, we propose to combine both within the same process. Indeed, we start with a random initial state and a small number of knots, then we make few iterations, then we choose a larger number of knots, interpolate the current spline at these new knots and then proceed in the optimization process. This strategy (Strategy I) has shown high accuracy and fast convergence.
\subsection{Strategy II: Automatic Knots Generation}
SPINS is based on interpolating splines at a few number of knots and evaluating a loss function at some collocation points. By default, the knots are considered uniformly distributed on the domain $[a,b]$. In some cases, considering the knots as the Chebeychev points provided faster convergence but this is not always true. From the possible shapes that can take the curve of a cubic spline, it seems that choosing the knots to be at the local extrema and inflection points of the unknown exact solution will provide better results. But since our algorithm starts the optimization process with a random initial state $\mathbf{y}$, we had the idea to start with uniformly distributed knots $(x_i)_i$ and a random initial state then make few iterations in order the current spline to take a stable shape that should be close but enough to the shape of the curve of the exact solution, then we detect its local extreme values and inflection points, and consider their locations as new knots. We apply this update every few iterations, then we continue the optimization until convergence. This strategy (Strategy II) has shown very good behaviour mainly when there are lots of variations (up and down) in the exact solution.\\

It is worth to mention that during the optimization process, it seems that SPINS learns the derivative functions of $S$ at early stages
before convergence. This observation gives credibility to the idea of updating the knots to be the locations where the derivatives of $S$ vanish. Here we talk about $S'$ and $S''$.\\

\begin{figure}[!htbp]
    \centering
    \includegraphics[width=0.32\textwidth]{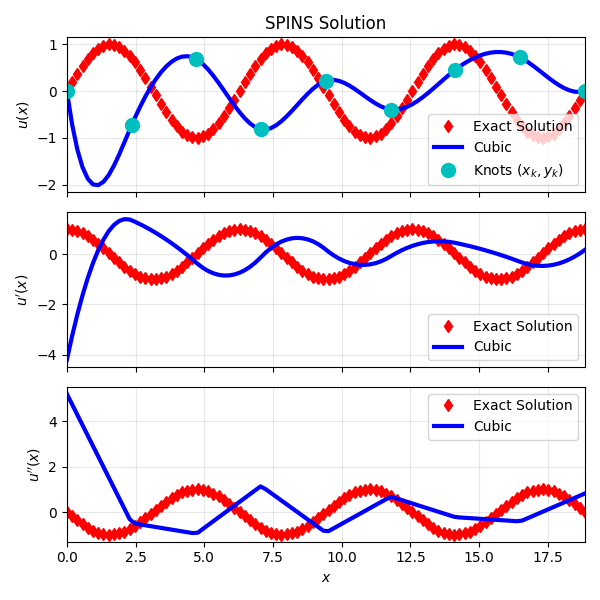}
    \includegraphics[width=0.32\textwidth]{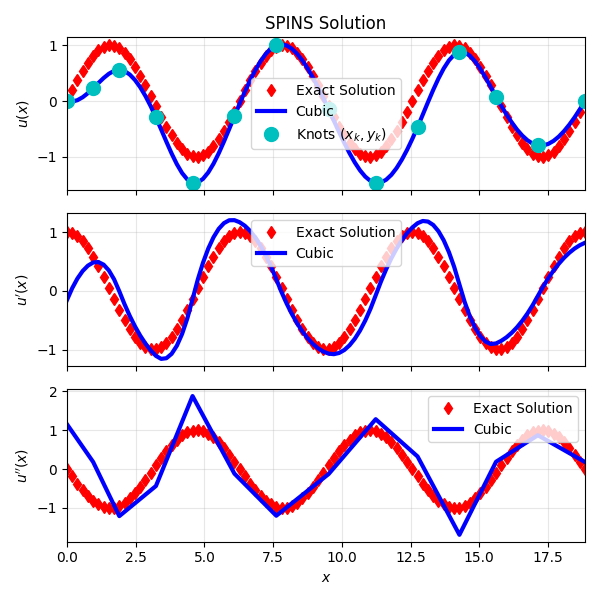}
    \includegraphics[width=0.32\textwidth]{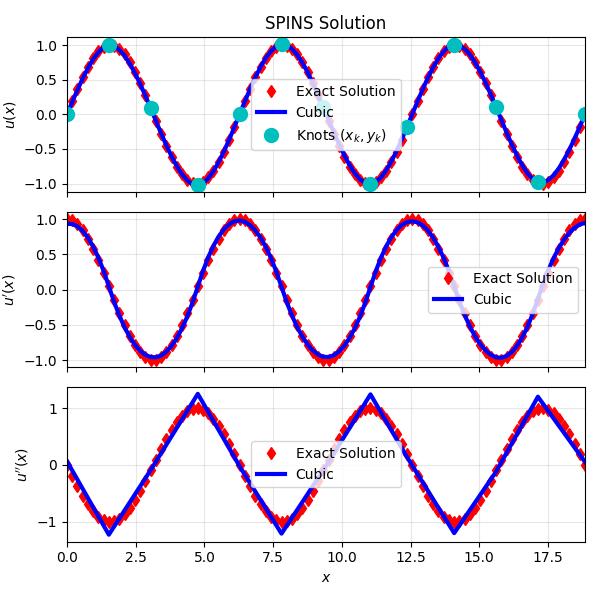}
    \caption{The Cubic SPINS and the exact solutions with their first and second order derivatives obtained using Strategy II. Left: the initialization of the optimization process with $9$ uniformly distributed knots and a random vector $\mathbf{y}$. Middle: the current spline solution after few iterations and detecting its local extrema and inflection points to be the new knots. Right: the final cubic SPINS solution.}
    \label{fig:StartegyII_Solution}
\end{figure}

For illustration of Strategy II, we consider BVP (\ref{BVP1}) but on the interval $[0,6\pi]$ with exact solution $u(x)=\sin(x)$. We solved this problem using cubic SPINS starting from $9$ uniformly distributed knots and random vector $\mathbf{y}$, then proceeding with the minimization for few iterations, detecting the local extrema and inflection points and consider them as new knots, updating the vector $\mathbf{y}$ from the latest spline, then completing the minimization process. The result is shown on Figure \ref{fig:StartegyII_Solution}.

When we do not have a sufficient number of new knots based on vanishing derivatives, this may happen depending on the shape of the exact solution $u$, then combining Strategy I and Strategy II works very well.
To illustrate the efficiency of these strategies combined together, we consider again BVP (\ref{BVP1}) with Dirichlet BCs on the interval $[0,3]$ with exact solution $u(x) =\dfrac{\sin(\pi x)}{1+x}$ and solve it using cubic SPINS.\\

\begin{figure}[!htbp]
    \centering
    \includegraphics[width=0.32\textwidth]{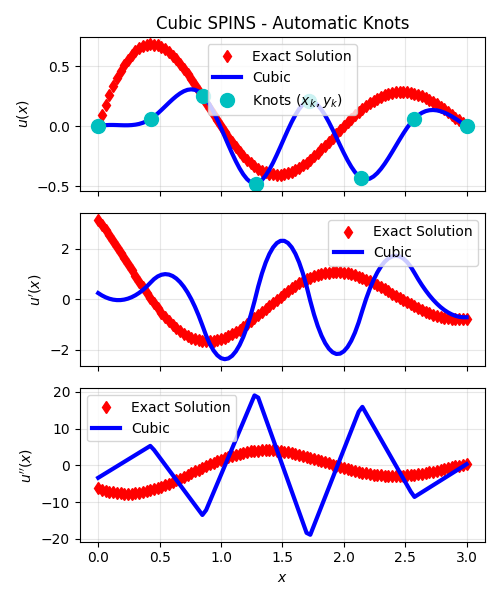}
    \includegraphics[width=0.32\textwidth]{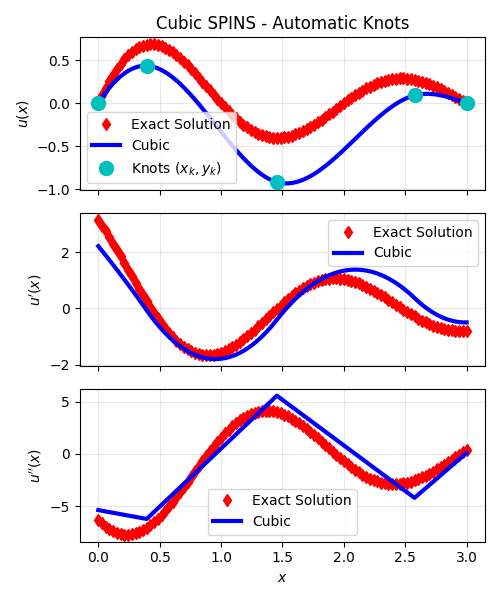}
    \includegraphics[width=0.32\textwidth]{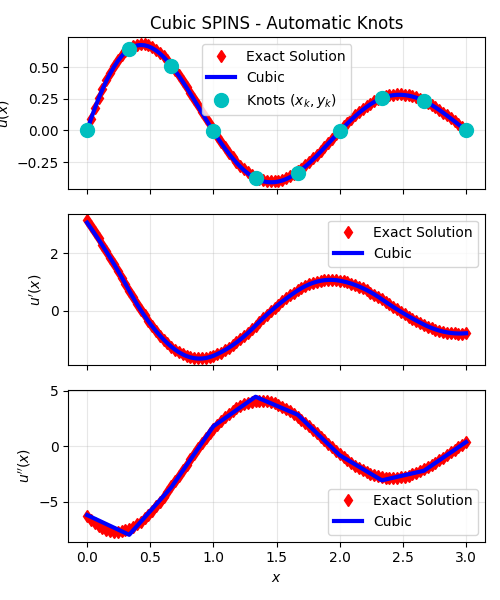}
    \caption{The Cubic SPINS and the exact solutions with their first and second order derivatives obtained by combining the strategies I and II. Left: the initialization of the optimization process with $8$ uniformly distributed knots and a random vector $\mathbf{y}$. Middle: the current spline solution after few iterations and detecting its local extrema to be the new knots. Right: the final SPINS solution with $10$ uniformly distributed knots.}
    \label{fig:example6_Solution}
\end{figure}

The result is shown on Figure \ref{fig:example6_Solution}, where we present the cubic spline obtained by interpolating $8$ uniformly distributed knots with a random initial state $\mathbf{y}$ using the BCs and the information of the derivatives at the boundaries extracted from the ODE. This cubic spline (left) is considered for the initialization of the optimization process. After making few iterations, the local extrema of the current spline are detected and their locations are considered to be the new knots. Again we make few iterations so that the spline takes a stable shape that should be close but not enough to the unknown exact solution (middle). Then we interpolate the current spline at $10$ uniformly distributed knots and we continue the optimization until convergence and get the final SPINS solution (right).

\paragraph{Remark} It is worth to mention that when computing the physics-informed residual at the collocation points, we obtain information where this residual is high, so we may refine the number of knots at these locations in order to compensate the high error.
\section{Conclusion and Future Work}
This paper introduces the Spline-based Physics-Informed Numerical Scheme (\textbf{SPINS}), a numerical framework designed to solve 1D Boundary Value Problems (BVPs). Traditional methods are physically robust but discrete, PINNs are continuous but boundary-weak. By utilizing high-order splines as basis functions, SPINS provides a continuous, high-order differentiable solution that strictly honors boundary conditions. This allows for analytical evaluation of the physics-informed residuals at any coordinate, ensuring that the physics is satisfied continuously. SPINS offers a robust, efficient, and mathematically transparent alternative to neural-network-based solvers. It bridges the gap between classical spline collocation and modern physics-informed optimization, providing a precise tool for solving 1D BVPs with minimal computational requirements.\\

While this work successfully demonstrates the efficacy of SPINS for 1D problems, several avenues remain for future exploration. Immediate next steps involve the extension of this framework to multi-dimensional domains (such as 2D and 3D PDEs) utilizing tensor-product splines, 2D B-spline surfaces or other smooth patched surfaces. Future extensions will also address evolution problems through a framework to strictly preserve temporal causality. Furthermore, a rigorous mathematical analysis is required to establish formal convergence proofs and to derive precise error bounds relative to the number of knots. Finally, investigating adaptive knot-placement strategies could further enhance the robustness of the solver for highly singular or boundary-layer dominated problems.

\paragraph{Declaration} This research did not receive any specific grant from funding agencies in the public, commercial, or not-for-profit sectors.\\
We have no conflicts of interest to disclose. 

\bibliographystyle{plain}
\bibliography{references} 

\end{document}